\documentclass{article}

\usepackage[active]{srcltx}

\usepackage{epsfig}
\usepackage{amsmath}
\usepackage{amsthm}
\usepackage{amssymb}
\usepackage{psfrag}
\usepackage[T1]{fontenc}
\usepackage{graphicx}
\usepackage[usenames,dvipsnames]{pstricks}
\usepackage[top=1in, bottom=1in, left=1in, right=1in]{geometry}
\usepackage{pst-grad} 
\usepackage{pst-plot} 

\graphicspath{{../figure/}}

\newtheorem{define}{Definition}

\newtheorem{prop}[define]{Proposition}

\newtheorem{theo}[define]{Theorem}
\newtheorem{cor}[define]{Corollary}
\newtheorem{lem}[define]{Lemma}
\newtheorem{prob}[define]{Problem}

\newcommand{\eps}{\varepsilon}
\newcommand{\Eps}{\mathcal{E}}

\newcommand{\w}{\omega}
\newcommand{\C}{\mathbf{C}}
\newcommand{\Z}{\mathbb{Z}}
\newcommand{\I}{\mathcal{I}}
\newcommand{\N}{\mathbb{N}}
\newcommand{\R}{\mathbb{R}}

\newenvironment{IEEEproof}{\begin{proof}}{\end{proof}}

\title{\LARGE \bf On the problem of reconstructing an unknown topology
		via locality properties of the Wiener filter}

\author{
	Donatello Materassi \quad
	Murti V.~Salapaka\\
		Department of Electrical and Computer Engineering,\\
		University of Minnesota,\\
		200 Union St SE, 55455, Minneapolis (MN) \\
		{\tt\small mater013@umn.edu}\qquad
		{\tt\small murtis@umn.edu}
		\medskip\\
}

\begin{document}

\maketitle

\begin{abstract}
The interest on networks of dynamical systems is increasing in the past years, especially because of their capability of modeling and describing a large variety of phenomena and behaviors.
Particular attention has been oriented towards the emergence of complicated phenomena from the interconnections of simple models.
We tackle, from a theoretical perspective, the problem of reconstructing the topology of an unknown network of linear dynamical systems where every node output is given by a scalar random process.
Existing approaches are bayesian or assume that the node inputs are accessible and manipulable. The approach we propose is completely blind, since we assume no a-priori knowledge of the network dynamics and we make use of second order statistics only. We also assume that it is only possible to observe the node outputs, while the inputs are not manipulable. The developed reconstruction technique is based on Wiener filtering and provides general theoretical guarantees for the detection of links in a network of dynamical systems. For a large class of networks, that we name ``self-kin'', sufficient conditions for the detection of the existing links are formulated. The necessity of those conditions is also discussed: indeed, it is shown that they do not hold only in specific pathological cases.
Thus, for any applications needs, the exact reconstruction of self-kin networks can be considered practically guaranteed.
For networks not belonging to this class those conditions are met by the smallest (in the sense of the least number of edges) self-kin network containing the actual one.
Hence, for general networks, the procedure identifies a self-kin network that is optimal in the sense of number of edges.
For networks not belonging to the self-kin class, conditions are met by the smallest (in the sense of the least number of edges) self-kin network containing the actual one.
\end{abstract}


\section{Introduction}
The interest on networks of dynamical systems is increasing in recent years, especially because of their capability of modeling and describing a large variety of phenomena and behaviors.
Particular attention has been oriented towards the emergence of complicated phenomena from interconnections of simple models (see, for example, \cite{CziBar99,LevRap00,FaxMur04,LiuYad08}).\\
Principal advantages provided by networked systems are three: a modular approach to design, the possibility of directly introducing redundancy, and the realization of distributed and parallel algorithms. All these advantages have led to increased utilization of realizationof  networked systems in the realization of many devices (see e.g. \cite{Gee05,KadBhu08}).\\
Interconnected systems are also successfully exploited to develop novel modeling approaches in many fields, such as
Economics (see e.g. \cite{NayRos07,ManSta00}),
Biology (see e.g. \cite{EisSpe98,RavSom02,DelNin08}),
Cognitive Sciences (see e.g. \cite{BroDin04}),
Ecology (see e.g. \cite{BunUrb00,UrbKei01})
and
Geology (see e.g. \cite{BaiMon06,MonBai05}),
especially when the investigated phenomena are characterized by spatial distributions where multivariate analysis is preferred \cite{InnMat08}.\\
While networks of dynamical systems are well studied and analyzed in physics \cite{BocLat06,GirNew02,NewGir04} and engineering \cite{ZhaLiu07,Olf07,SchRib08}, there are fewer results that address the problem of reconstructing an unknown dynamical network, since it poses formidable theoretical and practical challenges \cite{Kol09}.
However, unravelling the interconnectedness of a set of processes is of significant interest in many fields, and the necessity for general tools  is rapidly increasing (see \cite{Tim07}, \cite{BocIva07} and \cite{NapSau08} and the bibliography therein for recent results).
Existing results derive a network topology from sampled data (see e.g. \cite{ManSta00,Tim07,NapSau08,OzeUzu08}) or to determine the presence of substructures (see e.g. \cite{NewGir04,BocIva07}).
The Unweighted Pair Group Method with Arithmetic mean (UPGMA) \cite{MicSok57} is one of the first techniques proposed to reveal an unknown topology.
It has found widespread use  in the reconstruction of phylogenetic trees, and is widely employed in other areas such as communication systems and for resource allocation \cite{FreHar02}.
UPGMA identifies a tree topology relying on the observation of leaf nodes, theoretically guaranteeing a correct identification only on the strong assumption that an ultrametric is defined among the leaves.
Another well-known technique for the identification of a tree network is developed in \cite{ManSta00} for the analysis of a stock portfolio.
The authors identify a tree structure according to the following procedure: i) a metric based on the correlation index is defined among the nodes; ii) such a metric is employed to extract the Minimum Spanning Tree \cite{Die06} which forms the reconstructed topology.
Many improvements over \cite{ManSta00} have been devised especially using shrinking techniques to estimate the correlation matrix \cite{LedWol04,TumLil08} or refining its estimate via high frequency sampling \cite{TotKer09}. A reliability index for any link can also be defined using bootstrap techniques \cite{TumLil08}.
However, in \cite{InnMat08} a severe limit of this strategy is highlighted,
where it it is shown that, even though the actual network is a tree, the presence of dynamical connections or delays can lead to the identification of a wrong topology.
In \cite{MatInn09} a similar strategy, where the correlation metric is replaced by a metric based on the coherence function, is numerically shown to provide an exact reconstruction for tree topologies. Finally, in \cite{MatInn08} it is shown that a correct reconstruction can be guaranteed for any topology with no cycles.\\
An approach for the identification of more general topologies is developed in the area of Machine Learning for Bayesian Networks (BNs) \cite{GetFri02,FriKol03}.
The results presented in this paper present have strong connections with this area.
Indeed, the topology reconstruction of BNs is based on the fact that the conditional expectation operator to estimate one signal is determined only by the signals belonging to its Markov Blanket \cite{Pea88}. However, the most significant difference between BNs and the network models of this paper is that in BNs the network graph can not have loops.
The fact that the Wiener filter applied to the class of network models studied in this paper has properties similar to BNs even in presence of loops can not be considered trivial.
Moreover, we tackle the problem of the robustness of the reconstruction in terms of power spectral densities of the noise corrupting the signals.\\
In \cite{BocIva07} different techniques for quantifying and evaluating the modular structure of a network are compared and a new one is proposed trying to combine both the topological and dynamic information of the complex system. However, the network topology is only qualitatively estimated.
In \cite{Tim07} a method to identify a network of dynamical systems is described. However, primary assumptions of the technique are the possibility to manipulate the input of every single node and the possibility of conducting many experiments to detect the link connectivity.\\
In \cite{MarPel08a} an interesting and novel approach based on auto-regressive models and Granger-causality \cite{Gra69} is proposed for reconstructing a network of dynamical systems. 
This technique relies on multivariate identification procedure to detect the presence of a link, but still no theoretical sufficient or necessary conditions are derived to check the correctness of the results.\\
More recently, in \cite{NapSau08} and \cite{CanWak08} interesting equivalences between the identification of a dynamical network and a $l_0$ sparsification problem are highlighted, suggesting the difficulty of the reconstruction procedure \cite{CanTao05}.\\
Summarizing, to the best knowledge of the authors, apart from the results in \cite{MatInn08}, which are limited to tree topologies, no general theoretical guarantees about the correct reconstruction of network links are provided if there is no possibility of directly manipulating the system.\\
In this paper the problem of reconstructing a network of dynamical systems where every node represents an observable scalar signal and the dynamics is linear and represented by the connecting links is addressed.
The problem, when analyzed from a systems theory point of view, provides a method for correctly identifying a topology that belongs to the pre-specified class of self-kin networks. Moreover, if the network does not belong to such a class, conditions about the optimality of the identified topology according to a certain criterion is estabilished.
From this perspective, sufficient conditions for the exact reconstruction of a large class of networks, which we name self-kin, are derived. Examples of self-kin networks are given by (but not limited to) trees, and ring topologies \cite{JovBam05}. In the case the network is not self-kin, the reconstructed topology is guaranteed to be the smallest self-kin network containing the actual one.
The theory developed is not bayesian and relies directly on Wiener filtering theory.
Conditions derived for the detection of links are based on sparsity properties of the (non-causal) Wiener filter modeling the network.
Indeed, conditions under which the Wiener filter smoothing a signal of the network is ``local'' are derived.
From a different perspective, another important contribution of the paper is given by providing conditions for a local and distributed implementation of the Wiener filter. 

The paper is organized as follows.
In Section \ref{sec:intuition} examples are provided to provide the basic intuition behind central ideas;
in Section \ref{sec:preliminaries} definitions are provided based on standard notions of graph theory;
in Section \ref{sec:problem} the main problem is formulated;
in Section \ref{sec:non-causal} the main results are provided for non-causal Wiener filtering;
in Section \ref{sec:causal} the results are extended to causal Wiener filtering and Granger causality;
in Section \ref{sec:algorithm} the implementation of algorithms for the detection of network topologies are discussed for different scenarios;
in Section \ref{sec:corrupted} the robustness of the identification is addressed;
eventually, in Section \ref{sec:example} numerical simulations illustrating the effectiveness of the methodology are presented.

~\\~\\
\noindent{\bf Notation:}\\
The symbol $:=$ denotes a definition\\
$\|x\|$: $2$-norm of a vector $x$\\
$W^T$: the transpose of a matrix or vector $W$ \\
$W^*$: the conjugate transpose of a matrix or vector $W$ \\
$x_i$ or $\{x\}_i$: the $i$-th element of a vector $x$ \\
$W_{ji}$: the entry $(j,i)$ of a matrix $W$ \\
$W_{j*}$: j-th row of a matrix $W$\\ 
$W_{*i}$: i-th column of a matrix $W$\\
$x_V$: when $V=(v_1,...,v_n)$ is a $n$-tuple of natural numbers denotes the vector $(x_{v_1} ~...~x_{v_n})^T$\\
$|A|$: cardinality (number of elements) of a set $A$\\
$E[\cdot]$: mean operator; \\
$R_{XY}(\tau):= E[X(t)Y^T(t+\tau)]$: cross-covariance function of wide-sense stationary vector processes $X$ and $Y$;\\
$R_{X}(\tau):= R_{XX}(\tau)$: autocovariance; \\
$\mathcal{Z}(\cdot)$: Zeta-transform of a signal;\\
$\Phi_{XY}(z):=\mathcal{Z}(R_{XY}(\tau))$: cross-power spectral density;\\
$\Phi_{X}(z):=\Phi_{XX}(z)$: power spectral density;\\
with abuse of notation, $\Phi_{XY}(\w)=\Phi_{XY}(e^{i\w})$ for spectral densities and $W(\w)=W(e^{i\w})$ for transfer functions;\\
$b_{i}:$ $i$-th element of the canonical base of $\R^{n}$.\\

\begin{figure}[hbt!]
	\centering
	\includegraphics[width=0.6\columnwidth]{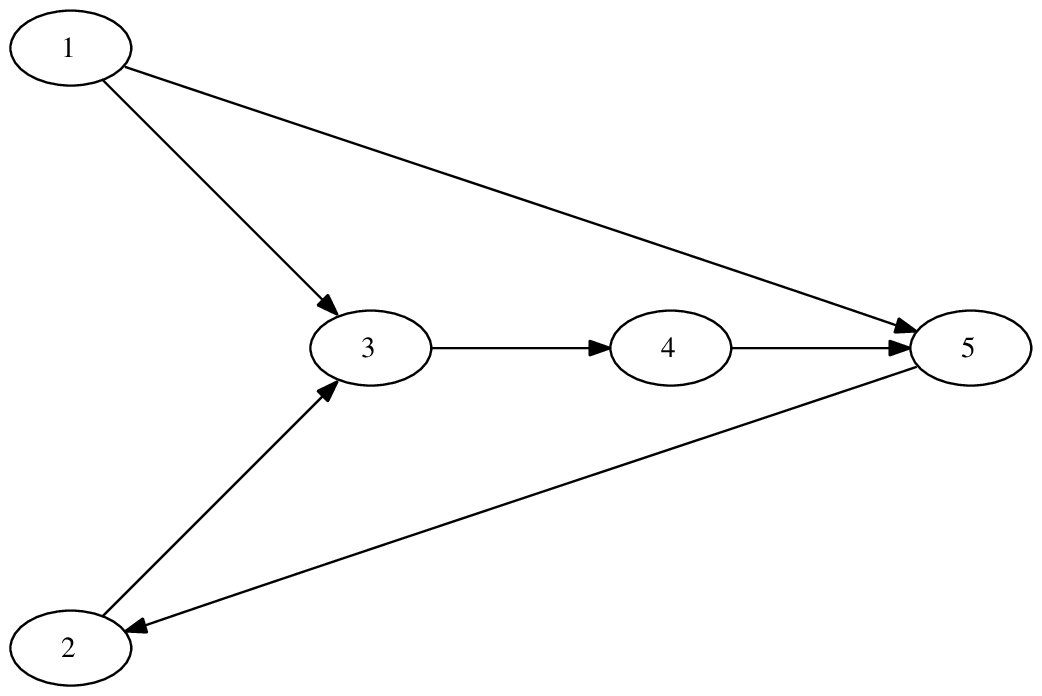}
	\caption{An example of the graphical representation of dynamical networks. Every node represents a stochastic signal and the edge from a node $N_i$ to a node $N_j$ represents the transfer function $H_{ji}(z)$. Every node signal is also implicitely affected by a process noise.\label{fig:example}}
\end{figure}

In this sections a representation of dynamical networks in terms of oriented graphs is presented.\\
In this representation, every node $N_j$ represents a scalar time-discrete wide-sense stationary stochastic process $x_j$, while every directed arc form a node $N_i$ to a node  $N_j$ represents a possibly non-causal transfer function $H_{ji}(z)\neq 0$.
The absence of such an arc implies that $H_{ji}(z)=0$.
Every node signal is also implicitly considered affected by an additive process noise $e_j$. Given the graphical representation, the dynamics of the network is described by
\begin{align}
	x_j = e_j + \sum_{i}H_{ji}(z)x_i.
\end{align}
For example, in Figure~\ref{fig:example} a network is represented, the dynamics of which corresponds to
\begin{align*}
	&x_1 = e_1\\
	&x_2 = e_2 + H_{25}(z)x_5\\
	&x_3 = e_3 + H_{31}(z)x_1+H_{32}(z)x_2\\
	&x_4 = e_4 + H_{43}(z)x_3\\
	&x_5 = e_5 + H_{51}(z)x_1+H_{54}(z)x_4.
\end{align*}

\section{Illustrative examples}\label{sec:intuition}
In this section special network configurations of dynamical systems are presented where the Wiener filter producing the estimate of a node signal $x_0$ from the other node signals $x_i$ has the characteristic of being sparse.
In most cases, the presence of a non-null entry in the Wiener filter corresponds to the presence of a direct link between the signal $x_i$ and the signal $x_0$.
In the following sections it will be shown that this result holds in general for the class of self-kin networks, a class of networks defined later.
Moreover if the network is not self-kin it will be shown that the presence of a non-null entry in the Wiener filter identifies the presence of a link in the ``smallest'' (in the sense of number of edges) self-kin network containing the original one.

\subsection{Wiener Filtering of a downstream signal}
\begin{figure}[hbt!]
	\centering
	\includegraphics[width=0.8\columnwidth]{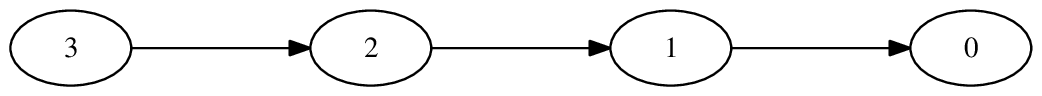}
	\caption{A cascade network. The Wiener filter estimating $x_0$ from the other signals makes use only of the ``parent'' signal $x_1$ which is directly linked to $x_0$.\label{fig:cascade dw}}
\end{figure}

Consider a network of four systems as represented in Figure~\ref{fig:cascade dw} where
\begin{align*}
	& x_{3}=e_{3}\\
	& x_{i-1}=e_{i-1}+H_{i-1,i}(z)x_{i}
\end{align*}
with $H_{i-1,i}(z)\neq 0$ three possibly non-causal SISO transfer functions for $i=1,2,3$ and with signals $e_i$ mutually uncorrelated for $i=0,1,2,3$.
Consider the Wiener filter that provides the estimate $\hat x_0$ of $x_0$ based on the other signals, $x_1$, $x_2$ and $x_3$. It can be shown that
\begin{align*}
	\hat x_{0}&=\left( W_{01}(z)~~ W_{02}(z)~~ W_{03}(z) \right)
		\left(\begin{array}{c}
			x_{1}\\ x_{2}\\ x_{3}
		\end{array}\right)=\\
		&=\left(H_{01}(z)~~ 0~~ 0 \right)
		\left(\begin{array}{c}
			x_{1}\\ x_{2}\\ x_{3}
		\end{array}\right)
\end{align*}
is the best estimate of $x_0$ in the least squares sense. Indeed,
\begin{align}
	E[(x_0 - \hat x_0)x_i]=E[e_0 x_i]=0
\end{align}
for $i=1,2,3$. Thus, from the Hilbert projection theorem \cite{Lue69}, 
\begin{align}
	W_0(z) := \left( W_{01}(z)~~ W_{02}(z)~~ W_{03}(z) \right)
\end{align}
is the Wiener filter providing the best estimate for $x_0$ \cite{Cai87}.
In this case, the only non-null entry of $W_0$ correctly detects that there is a connection between the node $0$ and the node $1$.

\subsection{Wiener Filtering of an upstream signal}
\begin{figure}[hbt!]
	\centering
	\includegraphics[width=0.8\columnwidth]{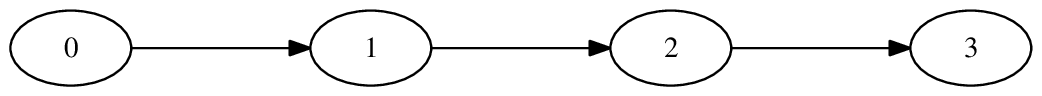}
	\caption{A cascade network. The Wiener filter estimating $x_0$ from the other signals makes use only of the ``child'' signal $x_1$ which is directly linked to $x_0$. \label{fig:cascade up}}
\end{figure}

Consider a network of four systems as represented in Figure~\ref{fig:cascade up} with
\begin{align}
	& x_{0}=e_{0}\\
	& x_{i+1}=H_{i+1,i}(z)x_{i}+e_{i+1}
\end{align}
where $H_{i+1,i}(z)$ are three possibly non-causal SISO transfer functions for $i=1,2,3$ and with signals $e_i$ mutually not correlated for $i=0,1,2,3$.
The Wiener filter that provides the estimate $\hat x_0$ of $x_0$ making use of the other signals $x_1$, $x_2$ and $x_3$ is such that
\begin{align*}
	\hat x_{0}&=\left( W_{01}(z)~~ W_{02}(z)~~ W_{03}(z) \right)
		\left(\begin{array}{c}
			x_{1}\\ x_{2}\\ x_{3}
		\end{array}\right)=\\
		&=\left( \frac{\Phi_{x_0 x_0}(z)H_{10}^*(z)}{|H_{10}(z)|^2\Phi_{x_0}(z)+\Phi_{e_1}}~~ 0~~ 0 \right)
		\left(\begin{array}{c}
			x_{1}\\ x_{2}\\ x_{3}
		\end{array}\right).
\end{align*}
Again in this case, there is one only non-null entry in the Wiener filter corresponding to the link between the node $0$ and the node $1$.

\subsection{Wiener Filtering of a loop}
\begin{figure}[hbt!]
	\centering
	\begin{tabular}{cc}
	\includegraphics[width=0.4\columnwidth]{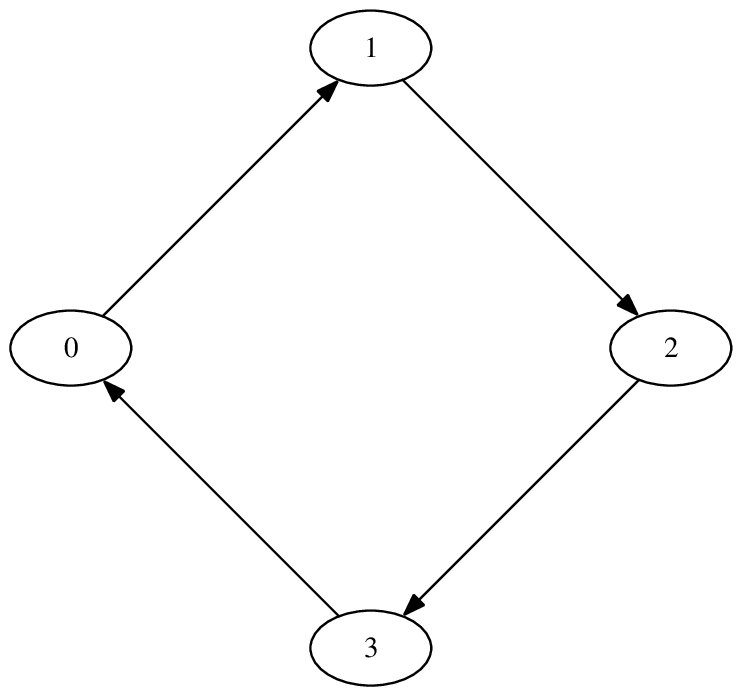} &
	\includegraphics[width=0.4\columnwidth]{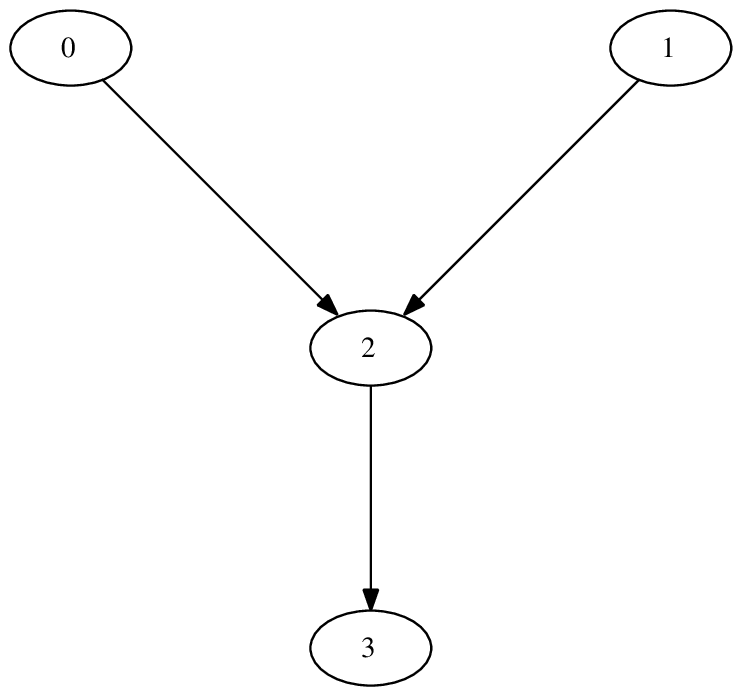}\\
	(a) & (b)
	\end{tabular}
	\caption{A loop network (a): the Wiener filter estimating $x_0$ from the other signals makes use only of the ``child'' signal $x_1$ and the parent signal $x_3$ which are directly linked to $x_0$\label{fig:loop} A more general network (b): the Wiener filter estimating $x_0$ from the other signals makes use of the signal $x_2$ which is directly linked to it, but also of the signal $x_1$ which is not.\label{fig:nonkin4}
	}
\end{figure}

In the two previous cases the network configurations did not involve loops and the resulting Wiener filter had the property of having non-null entries corresponding to the node signals immediately connected to the node of interest.
The identification of network topologies with loops is a challenging problem \cite{ManSta00, InnMat08}. Indeed, most techniques deal with networks with no cycles \cite{MatInn08}.
The presence of a loop leads to more complex relations between the node signals, especially in terms of the covariance function and (cross)-spectral densities.
Thus, it is interesting to note that the absence or presence of loops does not seem to affect the sparsity of the Wiener filter as shown in the following example \cite{MatInn08}.\\
Consider a network of four systems as represented in Figure~\ref{fig:loop}a where
\begin{align*}
	& x_{i}=H_{i,[i-1]_{mod~4}}(z)x_{[i-1]_{mod~4}}+e_{i}
\end{align*}
for $i=0,...,3$ and
\begin{align}
	& [n]_{mod~m}:=\min\{q| q=n+km\geq 0 \text{ and } k\in \Z \}.
\end{align}
Assume also that the signals $e_{i}$ are mutually not correlated.
It is possible to show that the best least squares estimate of $x_0$, based on $x_1$, $x_2$ and $x_3$, is given by
\begin{align*}
	&\hat x_{0}=\left( W_{01}(z)~~ W_{02}(z)~~ W_{03}(z) \right)
		\left(\begin{array}{c}
			x_{1}\\ x_{2}\\ x_{3}
		\end{array}\right)=\\
		&=\left( \frac{\Phi_{e_0}H_{10}(z)^{*}}{\Phi_{e_0}|H_{10}|^2+\Phi_{e_1}}~~ 0~~ (1-\frac{\Phi_{e_0}|H_{10}|^2}{\Phi_{e_0}|H_{10}|^2+\Phi_{e_1}})H_{03}(z) \right) \cdot \\
		&\qquad\cdot\left(\begin{array}{c}
			x_{1}\\ x_{2}\\ x_{3}
		\end{array}\right).
\end{align*}
It is to be noted that the signal $x_2$, which is not directly connected to $x_0$, is not used in the optimal estimate of $x_0$.
Thus, in all the examples presented so far the node signals actually employed by the Wiener Filter are the signals corresponding to nodes directly connected to the node signal being estimated.
This example, along with the previous ones, leads to the conjecture that the Wiener Filter can be used as a tool to identify which nodes are connected to a specific node of interest by checking wether the corresponding entries of the Wiener filter are null or not.
Unfortunately, the situation is more complex as shown in the following example.

\subsection{Wiener filter using a not directly connected signal}
Consider a network of four systems as represented in Figure~\ref{fig:nonkin4}b where
\begin{align}
	& x_{0}=e_{0}\\
	& x_{1}=e_{1}\\
	& x_{2}=H_{2,0}(z)x_{0}+H_{2,1}(z)x_{1}+e_{2}\\
	& x_{3}=H_{3,2}(z)x_{2}+e_{3}
\end{align}
and the noises $e_{j}$ are mutually not correlated.
In this example, the Wiener filter providing the estimate $\hat x_0$ of $x_{0}$ using the other signals of the network is given by
\begin{align*}
	\hat x_{0}=\left( W_{01}(z)~~ W_{02}(z)~~ W_{03}(z) \right)
		\left(\begin{array}{c}
			x_{1}\\ x_{2}\\ x_{3}
		\end{array}\right)=\\
		=\frac{\Phi_{e_0}(z)H_{20}^*(z)}{\Phi_{e_0}(z)|H_{20}(z)|^2+\Phi_{e_2}(z)}
		\left( -H_{21}(z)~~ 1~~ 0 \right)
		\left(\begin{array}{c}
			x_{1}\\ x_{2}\\ x_{3}
		\end{array}\right).
\end{align*}
In this case the Wiener filter estimating $x_0$ from the other signals makes use not only  of the signal $x_2$ which is directly linked to it, but also of the signal $x_1$ which is not.\\

After the analysis of these examples, it is natural to ask to what extent and under which assumptions it is possible to reconstruct the topology of a network of linear dynamical systems measuring the node signals and if the Wiener filter can be a useful tool for accomplishing this identification procedure.


\section{Preliminary definitions}\label{sec:preliminaries}
\begin{figure}[tb]
	\centering
	\begin{tabular}{cc}
		\includegraphics[width=0.45\columnwidth]{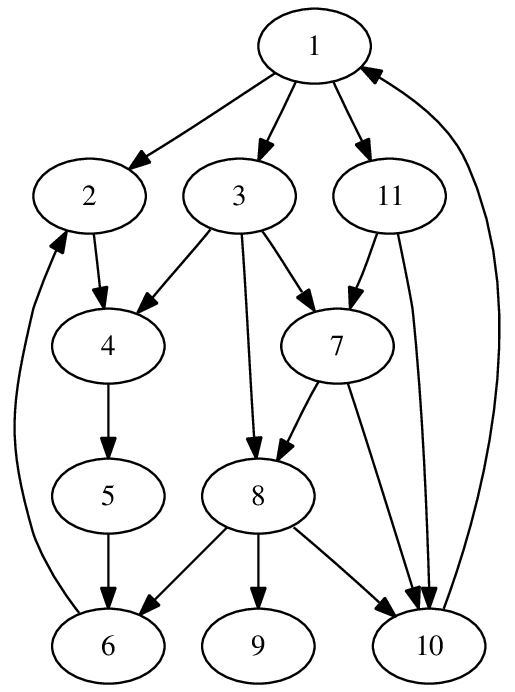} &
		\includegraphics[width=0.45\columnwidth]{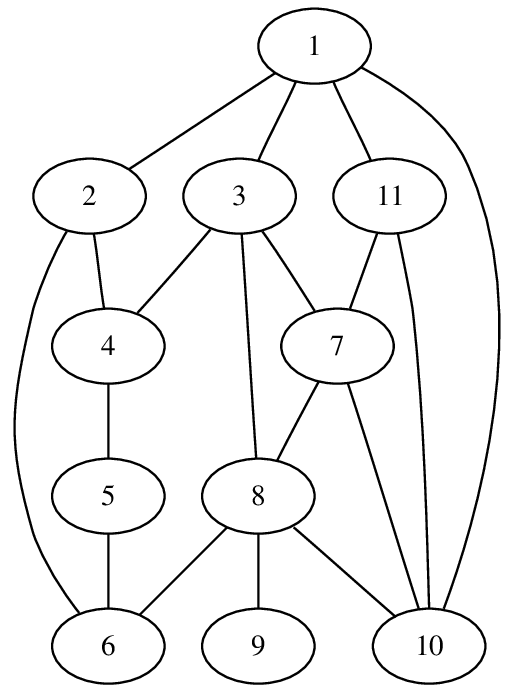}\\
		(a) & (b)
	\end{tabular}
	\caption{A directed graph (a) and its topology (b)\label{fig:topexample}.}
\end{figure}

In this section, basic notions of graph theory, which are functional to the following developments, will be recalled.
For an extensive overview see \cite{Die06}.
First, the standard definition of undirected and oriented graphs is provided.\\
\begin{define}[Directed and Undirected Graphs]
	An undirected graph $G$ is a pair $(V,A)$ where $V$ is a set of vertices or nodes and $A$ is a set of edges or arcs, which are unordered subsets of two distinct elements of $V$.\\
	A directed (or oriented) graph $G$ is a pair $(V,A)$ where $V$ is a set of vertices or nodes and  $A$ is a set of edges or arcs, which are ordered pairs of elements of $V$.\\
\end{define}
In the following, if not specified, oriented graphs are considered.\\
\begin{define}[Topology of a graph]
	Given an oriented graph $G=(V,A)$, its topology is defined as the undirected graph $G'=(V,A')$ such that $\{N_i,N_j\}\in A'$ if and only if $(N_i,N_j)\in A$ or $(N_j,N_i)\in A$, and $top(G):=G'$.\\
\end{define}
The topology of an oriented graph $G$ is the undirected graph $G'$ obtained by removing the orientation from any edge.
An example of a directed graph and its topology is represented in Figure~\ref{fig:topexample}.
~\\
\begin{define}[Children and Parents]
	Given a graph $G=(V,A)$ and a node $N_j\in~V$, we define the children of $N_j$ as
	$\mathcal{C}_G(N_j):=\{N_i| (N_j,N_i)\in A\}$ and the parents of $N_j$ as
	$\mathcal{P}_G(N_j):=\{N_i| (N_i,N_j)\in A\}$.\\
	Extending the notation, children and the parents of a set of nodes are denoted as follows
	\begin{align*}
		& \mathcal{C}_G(\{N_{j_1},...,N_{j_m}  \}):=
			\cup_{k=1}^m \mathcal{C}_G(N_{j_k})\\
		& \mathcal{P}_G(\{N_{j_1},...,N_{j_m}  \}):=
			\cup_{k=1}^m \mathcal{P}_G(N_{j_k}).
	\end{align*}
\end{define}
\begin{define}[Kins]
	Given an oriented graph $G=(V,A)$ and a node $N_j\in V$, kins of $N_j$ are defined as
	\begin{align*}
		\mathcal{K}_G(N_j):=\{N_i| N_i\neq N_j \text{ and }
					N_i \in \mathcal{C}_G(N_j)~\cup\\
			\qquad\qquad \cup~\mathcal{P}_G(N_j)~
						\cup~\mathcal{P}_G(\mathcal{C}_G(N_j)) \}.
	\end{align*}
	Kins of a set of nodes are defined in the following way
	\begin{align*}
		& \mathcal{K}_G(\{N_{j_1},...,N_{j_m}  \}):=
			\cup_{k=1}^m \mathcal{K}_G(N_{j_k}).
	\end{align*}
\end{define}
\begin{define}[Proper Parents and Children]
	Given an oriented graph $G=(V,A)$ and a node $N_j$, $N_i$ is a proper parent (child) of $N_j$ if it is a parent (child) of $N_j$ and
	$N_i\notin \mathcal{P}_G(\mathcal{C}_G(N_j))$.
	$N_i$ is a proper kin if it is a kin and
	$N_i\notin \mathcal{P}_G(N_j)\cup \mathcal{C}_G(N_j)$.\\
\end{define}
\begin{figure}[tb]
	\centering
	\begin{tabular}{cc}
		\includegraphics[width=0.5\columnwidth]{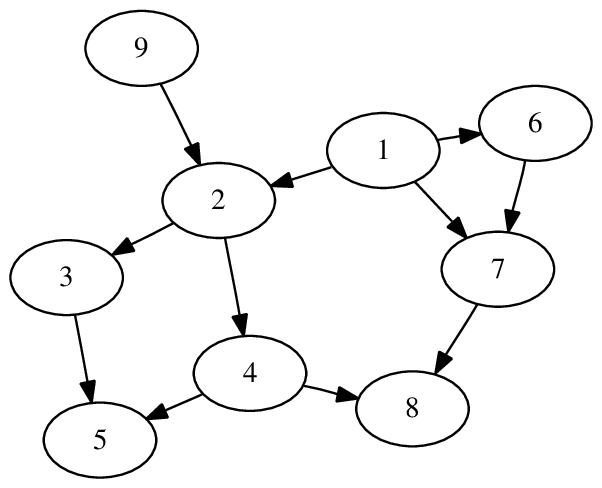} &
		\includegraphics[width=0.4\columnwidth]{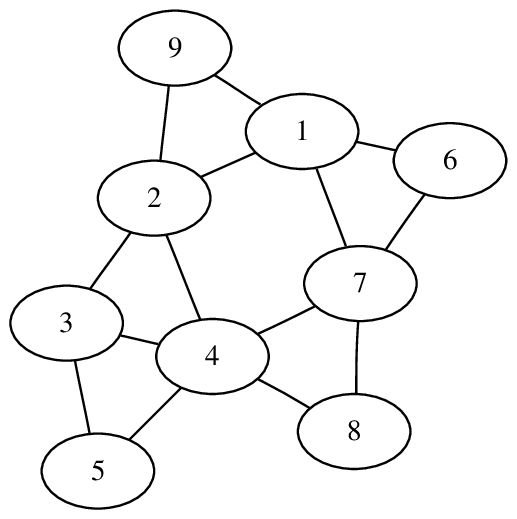}\\
		(a) & (b)
	\end{tabular}
	\caption{An oriented graph (a) and its kin topology (b)\label{fig:nonselfkin}.}
\end{figure}
Note that the kin relation is symmetric, in the sense that $N_i \in \mathcal{K_G}(N_j)$ if and only if $N_j \in \mathcal{K_G}(N_i)$.\\
\begin{define}[Kin-graph]
	Given an oriented graph $G=(V,A)$, its kin-graph is
	the undirected graph $\tilde G=(V,\tilde A)$ where
	\begin{align*}
		\tilde A := \{\{N_i, N_j\}| N_i \in \mathcal{K}_G(N_j)~ \text{for all } j \}.
	\end{align*}
	and it is denoted as $kin(G)=\tilde G$.\\
\end{define}
A directed graph and its kin-graph are represented in Figure~\ref{fig:nonselfkin}.
Note that the kin-graph of $G$ is an undirected graph. It could be defined as a directed graph, but, because of the symmetry of the kin relation, a directed graph contains exactly the same information. Moreover such a choice is motivated by the following definition
\begin{define}[Self-kin Graph]
	An oriented graph $G$ is self-kin if
	\begin{align*}
		top(G)=kin(G).
	\end{align*}
\end{define}

\begin{figure}[tb]
	\centering
	\includegraphics[width=0.75\columnwidth]{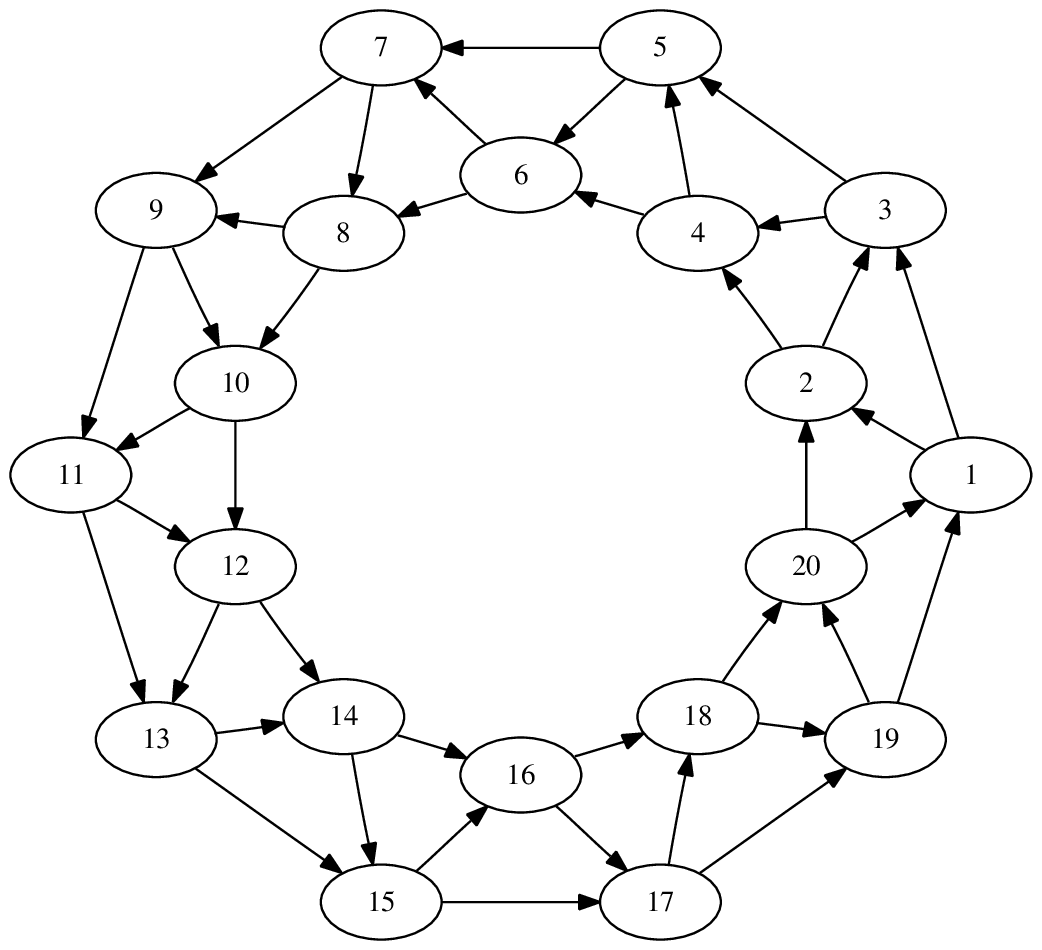}
	\caption{A ring network is always a self-kin network.\label{fig:selfkinring}}
\end{figure}

An example of self-kin network is given in Figure~\ref{fig:selfkinring}.
Many graphs are self-kin, such as oriented trees and rings \cite{Die06}.

\begin{define}
	Let $\Eps$ be a set containing time-discrete scalar, zero-mean, wide-sense, stationary random processes such that, for any $e_i,e_j\in\Eps$, the power spectral density $\Phi_{e_i e_j}(z)$ exists, is real rational with no poles on the unit circle and given by
	\begin{align}
		\Phi_{e_i e_j}(z)=\frac{A(z)}{B(z)},
	\end{align}
	where $A(z)$ and $B(z)$ are polynomials with real coefficients such that $B(z)\neq 0$ for any $z\in \C$, with $|z|=1$. 
	Then, $\Eps$ is a set of rationally related random processes.\\
\end{define}

\begin{define}
	The set $~{}^0\mathcal{F}$ is defined as the set of real-rational SISO transfer function with no poles on the unit circle $\{z\in\C|~~|z|=1\}$.\\
\end{define}

\begin{define}
	Given a SISO transfer function $H(z)\in~{}^0\mathcal{F}$, represented as
	\begin{align}
		H(z)=\sum_{k=-\infty}^{\infty}h_k z^{-k},
	\end{align}
	we define the causal truncation operator as
	\begin{align}
		\{H(z)\}_{C}:=\sum_{k=0}^{\infty}h_k z^{-k}.
	\end{align}
\end{define}

\begin{lem}
	For every $H(z)\in~{}^0\mathcal{F}$, we have that $\{H(z)\}_{C}\in~{}^0\mathcal{F}$.
\end{lem}

\begin{define}
	The set $~{}^0\mathcal{F}^{+}$ is defined as the set of real-rational SISO transfer functions in $~{}^0\mathcal{F}$ such that
	\begin{align}
		\{H(z)\}_{C}=H(z).
	\end{align}
\end{define}

\begin{define}
	Let $\Eps$ be a set of rationally related random processes.
	The set ${}^{0}\mathcal{F}\Eps$ is defined as
	\begin{align*}
		{}^0\mathcal{F}\Eps:=
		\left\{x=\sum_{k=1}^{m} H_k(z)e_k~|~e_k\in\Eps, H_k(z)\in {}^0\mathcal{F}, m\in\N \right\}.
	\end{align*}
\end{define}
\begin{lem}
	The set ${}^0\mathcal{F}\Eps$ is a vector space with the field of real numbers. Let
\begin{align*}
	<x_1,x_2>:=R_{x_1 x_2}(0)=\int_{-\pi}^{\pi}\Phi_{x_1 x_2}(\w),
\end{align*}
which defines an inner product on ${}^0\mathcal{F}\Eps$ with the assumption that two processes $x_1$ and $x_2$ are considered identical if $x_1(t)=x_2(t)$, almost always for any $t$.
\end{lem}
\begin{IEEEproof}
      See supplemental material.
\end{IEEEproof}

For any $x\in {}^0\mathcal{F}\Eps$, the norm induced by the inner product is defined as $\|x\|:=\sqrt{<x,x>}$.
\begin{define}
	For a finite number of elements $x_1,...,x_m\in{}^0\mathcal{F}\Eps$, tf-span is defined as
	\begin{align*}
		&\text{tf-span}\{x_1,...,x_m\}:=\\
		&:=\left\{x=\sum_{i=1}^m \alpha_{i}(z)x_i
		~\lvert~\alpha_{i}(z)\in{}^0\mathcal{F}
		\right\}.
	\end{align*}
\end{define}
\begin{lem}
	The tf-span operator defines a subspace of ${}^0\mathcal{F}\Eps$.\\
\end{lem}
\begin{IEEEproof}
      The proof is left to the reader.
\end{IEEEproof}

\begin{define}
	For a finite number of elements $x_1,...,x_m\in{}^0\mathcal{F}\Eps$, c-tf-span is defined as
	\begin{align*}
		&\text{c-tf-span}\{x_1,...,x_m\}:=\\
		&:=\left\{x=\sum_{i=1}^m \alpha_{i}(z)x_i
		~\lvert~\alpha_{i}(z)\in{}^0\mathcal{F}^{+}
		\right\}.
	\end{align*}
\end{define}
\begin{lem}
	The c-tf-span operator defines a subspace of ${}^0\mathcal{F}\Eps$.\\
\end{lem}
\begin{IEEEproof}
      The proof is left to the reader.
\end{IEEEproof}

The following definition provides a class of models for a network of dynamical systems.
\begin{define}
	Let $I=(1,...,n)^T$ and consider the triplet $\mathcal{G}=(G,\mathcal{H},e_I)$ where $G=(V,A)$ is a graph with $n$ vertexes $\{N_j\}_{j=1,...,n}$; $\mathcal{H}:A\rightarrow~ {}^0\mathcal{F}$ is a function associating a, possibly non-causal, transfer function $H_{ji}(z)\in {}^0\mathcal{F}$, with no poles on the unit circle, to any edge $(N_i, N_j)\in A$; and $\Eps:=e_I=(e_1,...,e_n)^T$ is an ordered $n$-tuple of rationally related random processes.
	$\mathcal{G}$ is a Linear Dynamic Graph (LDG) if the following conditions are satisfied
	\begin{itemize}
		\item	$(N_i,N_j)\in A$ implies $(N_j,N_i)\notin A$
		\item	$\Phi_{e_i e_j}(\w)=0$ for $i\neq j$ and for all $\w\in \R$ 
	\end{itemize}
	The output $x_j$ of the dynamic of a LDG is defined as
	\begin{align*}
		x_j=e_j+\sum_{i}H_{ji}(z)x_i.
	\end{align*}
	for any $j=1,...,n$.
	In a matrix form it is possible to write the network dynamics of a LDG as
	\begin{align}
		x_{I}=e_{I}+H(z)x_{I}
	\end{align}
	where the entry $(j,i)$ of the transfer matrix $H(z)$ is given by the transfer function $H_{ji}(z)$.
	The LDG is {\it well-posed} if the matrix $(\I-H(z))$ is invertible, and both $(\I-H(z))$ and $(\I-H(z)) ^{-1}$ have full normal rank and no poles on the unit circle.
	This implies that there is a transfer matrix $T(z)=[\I-H(z)]^{-1}$ such that $x_{I}=T(z)e_{I}$. Interpreting $T(z)$ as a possibly non-causal, linear transformation, it follows that $x_{i}\in {}^0\mathcal{F}\Eps$ for any $i=1,...,n$.\\
	The LDG is {\it causal} if each entry of $H(z)$ and each entry of $T(z)$ is in ${}^0\mathcal{F}^{+}$.
\end{define}
A LDG is a complex interconnection of linear transfer functions $H_{ji}(z)$ connected according to a graph $G$ and forced by stationary additive mutually uncorrelated noise.
The following definitions will be useful for determining sufficient conditions for detection of links in a network.
\begin{define}
	A LDG $\mathcal{G}=(G,\mathcal{H},e_I)$ is topologically detectable if
	$\Phi_{e_i}(\w)>0$ for any $\w\in[-\pi,\pi]$ and for any $i=1,...,n$.\\
\end{define}

\section{Problem Formulation}\label{sec:problem}
\begin{prob}
	Consider a well-posed LDG $\mathcal{G}=(G,\mathcal{H},\mathcal{E})$ where the graph $G$ is unknown.
	Given the Power (Cross-)~Spectral Densities of $\{x_j\}_{i=1,\ldots,n}$,
	reconstruct the unknown topology of the graph $G$ of the LDG $\mathcal{G}$.
\end{prob}

\section{Sparsity of the non-causal Wiener Filter}\label{sec:non-causal}

First, a lemma is provided that guarantees that any element in $\text{tf-span}\{x_i\}_{i=1,...,n}$ admits a unique representation if the cross-spectral density matrix of its generating processes has full normal rank.
\begin{lem}\label{lem:unique representation}
	Let $q$ and $x_1,...,x_n$ be processes in the space $^0\mathcal{F}\Eps$.
	Define $I:=(1,....,n)$, and
	suppose that $x\in\text{tf-span}\{x_i\}_{i=1,...,n}$ and that $\Phi_{x_{I}x_{I}}(\w)>0$ almost for any $\w\in~[-\pi, \pi]$. Then there exists a unique transfer matrix $H(z)$ such that $q=H(z)x_{I}$.
\end{lem}
\begin{IEEEproof}
	Note that if $H(z)$ is such that $H(z)x_I=0$, then
	\begin{align}
		\Phi_{qq}(\w)=0=H(\w)\Phi_{x_{I}x_{I}}(\w)H^{*}(\w).
	\end{align}
	Since $\Phi_{x_{I}x_{I}}(\w)>0$ almost everywhere, we have $H(\w)=0$ almost everywhere which implies that $H(z)=0$.
	Now, by contradiction assume that $q=H_{1}(z)x_{I}=H_{2}(z)x_{I}$, with $H_1(z)\neq H_2(z)$. Then $0=[H_2(\w)-H_1(\w)]\Phi_{x_I x_I}(\w)[H_2(\w)-H_1(\w)]^*$ implying that $H_1(z)=H_2(z)$.
\end{IEEEproof}
A specific formulation of the non-causal Wiener filter is introduced for the defined spaces.
\begin{prop}\label{eq: my wiener}
	Let $x$ and $x_1,...,x_n$ be processes in the space $^0\mathcal{F}\Eps$.
	Define $I:=(1,....,n)$ and $X:=\text{tf-span}\{x_1,...,x_n\}$. Consider the problem
	\begin{align}\label{eq: cost general wiener}
		\inf_{q \in X} \|x-q\|^2.
	\end{align}
	If $\Phi_{x_I}(\w)>0$, for $\w\in[-\pi,\pi]$, the solution $\hat x\in X$ exists, is unique and has necessarily the form $\hat x=W(z)x_I$ where
	\begin{align*}
		W(z)=\Phi_{x x_{I}}(z)\Phi_{x_{I}}(z)^{-1}.
	\end{align*}
	Moreover, it holds that $\hat x$ is the only element in $X$ satisfying, for any $q\in X$,
	\begin{align}\label{eq:perp hilbert projection thm}
		<x-\hat x,q>=0.
	\end{align}
\end{prop}
\begin{IEEEproof}
	Observe that, since $q\in X$, the cost function satisfies
	\begin{align*}
		&\|x-W(z)x_{I}\|^2=\\
		&\quad=\int_{-\pi}^{\pi}\Phi_{xx}(\w)+W(\w)\Phi_{x_{I}x_{I}}(\w)W^*(\w)+\\
		&\qquad-\Phi_{x_{I}x}(\w)W^{*}(\w)-W(\w)\Phi_{xx_{I}}(\w).
	\end{align*}
	The integral is minimized by minimizing the integrand for all $\w\in [-\pi,\pi]$.
	It is straightforward to find that the minimum is achieved for
	\begin{align*}
		W(\w)=\Phi_{x x_{I}}(\w)\Phi_{x_{I} x_{I}}(\w)^{-1}.
	\end{align*}
	Defining the filter $W(z)=\Phi_{x x_{I}}(z)\Phi_{x_{I} x_{I}}(z)^{-1}$ a real-rational transfer matrix is obtained with no poles on the unit circle that has the specified frequency response. Thus $\hat x=W(z)x_I\in X$ minimizes the cost (\ref{eq: cost general wiener}).
	As a consequence of the Hilbert projection theorem (for pre-Hilbert spaces) Equation (\ref{eq:perp hilbert projection thm}) is satisfied for $\hat x$ if and only if it is the unique element of the subspace $X$ minimizing (\ref{eq: cost general wiener}) \cite{Lue69}.
	If $\Phi_{x_{I} x_{I}}(z)>0$, the uniqueness of $W(z)$ follows from Lemma \ref{lem:unique representation}.
\end{IEEEproof}
In the following definition a notion of conditional non-causal Wiener-uncorrelation is given.
\begin{define}
	Let $q$, $x_1,...,x_n$ be processes in the space $^0\mathcal{F}\Eps$.
	Define $I:=(1,....,n)$ and $X:=\text{tf-span}\{x_1,...,x_{n}\}$.
	For any $i\in\{1,...,n\}$, the process $x_i$ is non-causally Wiener-uncorrelated with $q$ given the processes $\{x_{k}\}_{k\neq i}$ if the $i$-th entry of the Wiener filter to estimate $q$ from $x_I$ is zero, that is
	\begin{align}
		\Phi_{q x_I}\Phi_{x_I x_I}^{-1}b_i=0.
	\end{align}
\end{define}
The following lemma provides an immediate relationship between non-causal Wiener-uncorrelation and the inverse of the cross-spectral density matrix.
This result presents strong similarities with the property of the inverse of the covariance matrix for jointly Gaussian random-variables. Indeed, it is well-known that the entry $(i,j)$ of inverse of the covariance matrix of $n$ random variables $x_1,...,x_n$ is zero if and only if $x_i$ and $x_j$ are conditionally independent given the othe variables.
\begin{lem}\label{lem:topological filter}
	Let $x_1,...,x_n$ be processes in the space $^0\mathcal{F}\Eps$.
	Define $I:=(1,....,n)$. Assume that $\Phi_{x_I}$ has full normal rank.
	The process $x_i$ is non-causally Wiener-uncorrelated with $x_j$ given the processes $\{x_k\}_{k\neq i,j}$, if and only if
	the entry $(i,j)$, or equivalently the entry $(j,i)$, of $\Phi_{x_I}^{-1}(z)$ is zero, that is, for $i\neq j$,
	\begin{align}
		b_{j}^{T}\Phi_{x_I}^{-1}b_{i}=b_{i}^{T}\Phi_{x_I}^{-*}b_{j}=0.
	\end{align}
\end{lem}
\begin{IEEEproof}
	Without any loss of generality, let $j=n$.
	Define the vector $x_{\overline{n}}=(x_1,...,x_{n-1})^T$ and determine the non-causal Wiener filter $W_{n\overline{n}}$ estimating $x_n$ from $x_{\overline{n}}$.
	\begin{align}
		x_{n}=\eps_n +W_{n\overline{n}}(z)x_{\overline{n}}
	\end{align}
	where the error $\eps_n$ has the property that $\Phi_{\eps_n x_{\overline{n}}}(z)=0$. Define $r:=(x_{\overline{n}^T}, \eps_n )$ and observe that
	\begin{align}
		r=\left(\begin{array}{cc}
			\I & 0\\ -W_{n\overline{n}}(z) & 1
		\end{array}\right)x; \qquad
		x_{I}=\left(\begin{array}{cc}
			\I & 0\\ W_{n\overline{n}}(z) & 1
		\end{array}\right)r.
	\end{align}
	Compute
	\begin{align}
		&\Phi_{x_I}^{-1}=
		\left(\begin{array}{cc}
			\I & W_{n\overline{n}}(z)^{*}\\ 0 & 1
		\end{array}\right)
		\left(\begin{array}{cc}
			\Phi_{x_{\overline{n}}}^{-1} & 0\\ 0 & \Phi_{\eps_{n}}^{-1}
		\end{array}\right)
		\left(\begin{array}{cc}
			\I & 0\\ W_{n\overline{n}}(z) & 1
		\end{array}\right)=\\
		&=
		\left(\begin{array}{cc}
			\Phi_{x_{\overline{n}}}+
			\frac{W_{n\overline{n}}^*W_{n\overline{n}}}{\Phi_{\eps_n}}
			& W_{n\overline{n}}^* \Phi_{\eps_n}^{-1}\\
			\Phi_{\eps_n}^{-1} W_{n\overline{n}} & \Phi_{\eps_n}^{-1}
		\end{array}\right).
	\end{align}
	The assertion is proven by premultiplying by $b_n^T$ and post-multiplying by $b_{i}$
\end{IEEEproof}
The following theorem provides a sufficient condition to determine if two nodes in a LDG are kins.
\begin{theo}\label{thm:kinship non-null}\label{thm:sparse wiener}
	Consider a well-posed and topologically detectable LDG.
	Let $x_1,...,x_n\in {}^0\mathcal{F}\Eps$ be the signals associated with the $n$ nodes of its graph.
	Define $X_j=\text{tf-span}\{x_i\}_{i\neq j}$.
	Consider the problem of approximating the signal $x_j$ with an element $\hat x_j\in X_j$, as defined below
	\begin{align}\label{eq:min problem}
		\inf_{\hat x_j\in X_j}
		\left\|x_j-\hat x_j \right\|^2.
	\end{align}
	Then the optimal solution $\hat x_j$ exists, is unique and
	\begin{align}
		\hat x_j =\sum_{i\neq j} W_{ji}(z)x_i
	\end{align}
	where $W_{ji}(z)\neq 0$ implies $(N_i, N_j) \in kin(G)$.
\end{theo}
\begin{IEEEproof}
	The LDG dynamics is given by $x=(\I-H(z))^{-1}e$ implying that $\Phi_{xx}$ can also be represented as
	\begin{align}
		\Phi_{xx}^{-1}=(\I-H)^{*}\Phi_{e}^{-1}(\I-H).
	\end{align}
	Consider the $j$-th line of $\Phi_{xx}^{-1}$. We have
	\begin{align}
		b_{j}^{T} \Phi_{xx}^{-1}=(b_{j}^{T}-H_{*j}^*)\Phi_{e}^{-1}(\I-H).
	\end{align}
	The row vector $(b_{j}^{T}-H_{*j}^*)$ has zero entries corresponding any index $k$ such that $k\neq j$ that is not a parent of $j$.
	Since $\Phi_{e}$ is diagonal the $i$-th column of $\Phi_{e}^{-1}(\I-H)$ has zero entries for any $k\neq i$ and is not a parent of $i$. Given $i\neq j$, if $i$ is not a parent of $j$ and $i$ is not a child of $j$ and $i$ and $j$ have no common children (they are not coparents), it follows that the entry $(j,i)$ of $\Phi_{xx}^{-1}$ is not zero. Using Lemma (\ref{lem:topological filter}) the assertion is proven.
\end{IEEEproof}
The following result provides a sufficient condition for the reconstruction of a link in a LDG.
\begin{cor}
	Consider a LDG $\mathcal{G}=(G,\mathcal{H},e_I)$ with $G=(V,A)$
	Assume $V=\{N_i\}_{i=1,...,n}$ and let $X:=\{x_i\}_{i=1,...,n}$ be the set of processes $x_i$ corresponding to the nodes $N_i$.
	Let $W_{ji}(z)$ be entry of the non-causal Wiener filter estimating $x_j$ from $\{x_k\}_{k\neq j}$ corresponding to the process $x_i$.
	If $\mathcal{G}$ is self-kin, then $W_{ji}(z)\neq 0$ implies $(N_j, N_i) \in top(G)$.
\end{cor}
\begin{IEEEproof}
	Since $\mathcal{G}$ is self-kin,
	$\mathcal{P}_G(N_j)\cup \mathcal{C}_G(N_j) \cup \mathcal{P}_G(\mathcal{C}_G(N_j))=
		\mathcal{C}_G(N_j) \cup \mathcal{P}_G((N_j))$.
	Thus, from the previous theorem the assertion follows immediately.
\end{IEEEproof}

The following lemma is a key result to explicitly determine the expression of the Wiener filter for a LDG in the non-causal and in the causal scenarios.
\begin{lem}\label{lem:smart span}
	Consider a well-posed LDG $\mathcal{G}=(G,\mathcal{H},e_I)$.
	Assume $V=\{N_i\}_{i=1,...,n}$ and let $X:=\{x_i\}_{i=1,...,n}$ be the set of processes $x_i$ corresponding to the nodes $N_i$.
	Fix $j\in\{1,...,n\}$ and define the set
	\begin{align*}
		& C:=\{c| N_c\in \mathcal{C}_G(N_j)\}=\{c_1, ... c_{n_c}\}
	\end{align*}
	containing the indexes of the $n_c$ children of $N_j$.
	Then, for $i \neq j$,
	\begin{align*}
		x_i \in \text{tf-span}
		\left\{
			\left\{
				\bigcup_{k \in C } (e_k+H_{kj}(z)e_j)
			\right\}
			\cup
			\left\{
				\bigcup_{k \notin C \cup \{j\}} \{e_k\}
			\right\}
		\right\}.
	\end{align*}
	Furthermore, if $\mathcal{G}$ is causal,
	\begin{align*}
		x_i \in \text{c-tf-span}
		\left\{
			\left\{
				\bigcup_{k \in C } (e_k+H_{kj}(z)e_j)
			\right\}
			\cup
			\left\{
				\bigcup_{k \notin C \cup \{j\}} \{e_k\}
			\right\}
		\right\}.
	\end{align*}
\end{lem}
\begin{IEEEproof}
	Define
	\begin{align}
		&\eps_{j}:=0 \nonumber\\
		&\eps_{k}:=e_k+H_{kj}(z)e_j	&\text{if } k \in C\nonumber\nonumber\\
		&\eps_{k}:=e_k		&\text{if }k \notin \{C \} \cup \{j\}\nonumber\\
		&\xi_{k}:=\sum H_{ki}(z)x_{i}	&\text{if } k=j&\nonumber\\
		&\xi_{k}:=x_k		&\text{if }k \neq j \label{eq:xi when not i}
	\end{align}
and, by inspection, observe that
\begin{align*}
	[\I-H(z)]\xi_{I}=\eps_{I}.
\end{align*}
	Since $\mathcal{G}$ is well posed, $[\I-H(z)]$ is invertible implying that the signals $\{\xi_i\}_{i=1,...,n}$ are a linear transformation of the signals $\{\eps_i\}_{i=1,...,n}$. For $i\neq j$, we have
\begin{align*}
	x_i = \xi_i \in \text{tf-span}\{\eps_k\}_{k=1,...,n}=\text{tf-span}\{\eps_k\}_{k\neq j}
\end{align*}
	where the first equality follows from (\ref{eq:xi when not i}) and  last equality follows from the fact that $\eps_{j}=0$.
	The causality of $\mathcal{G}$ also implies that
	\begin{align*}
		x_i =\text{c-tf-span}\{\eps_k\}_{k\neq j}.
	\end{align*}
	This proves the assertion.
\end{IEEEproof}
In the case of a well-posed and topologically detectable LDG, the following theorem
completes Theorem~\ref{thm:kinship non-null}
providing an explicit expression for the non-causal Wiener filter $W$ since the solution of the minimization problem (\ref{eq:min problem}) is known to be unique.
\begin{theo}\label{thm:kinship explicit}
	Consider a well posed and topologicaly detectable LDG $\mathcal{G}=(G,\mathcal{H},e_I)$ with $G=(V,A)$, $V=\{N_1,...,N_n\}$.
	Define $I=(1,...,n)^T$ and $I_j=(1,...,j-1,0,j+1,...,n)$.
	Consider an auxiliary process $e_0$ such that $\Phi_{e_0 e_I}=0$ and $\Phi_{e_0}=1$.
	Define,
	\begin{align}\label{eq: child filter}
		\hat C_{ji}(z)=
			\{\Phi_{e_j}H^{*}_{*j}(z)
			[\Phi_{e_j}(z)H_{*j}(z)H^{*}_{*j}(z)+\Phi_{e_{I_j}}(z)]^{-1}\}_i
	\end{align}
	\begin{align}\label{eq: parent filter}
		\hat P_{ji}(z)=
			\left(
				1-\hat C_{j*}(z)H_{*j}(z)
			\right)H_{ji}(z)
	\end{align}
	\begin{align}\label{eq: kin filter}
		\hat K_{ji}(z)=
			\left\{\begin{array}{ll}
				0 & \text{  if } i=j\\
				-\hat C_{j*}(z)H_{*i}(z) & \text{otherwise}.
			\end{array}\right.
	\end{align}
	These definitions do not depend on the choice on $e_0$.
	Consider the estimate $\hat x_j$ for $x_j$ provided by the non-causal Wiener filter using the other signals $\{x_i\}_{i\neq j}$
	\begin{align}
		\hat x_{j}=\sum_{i\neq j}W_{ji}x_i.
	\end{align}
	Then, we have that
	\begin{align}
		W_{ji}(z)= \hat C_{ji}(z)+ \hat P_{ji}(z)+\hat K_{ji}(z).
	\end{align}
\end{theo}
\begin{IEEEproof}
	Define the following set of indexes
	\begin{align*}
		& C:=\{c| N_c\in \mathcal{C}_G(N_j)\}=\{c_1, ... c_{n_c}\}\\
		& P_l:=\{p| N_p\in \mathcal{P}_G(N_l)\}/\{j\}=\{p_{l1}, ... p_{ln_{pl}}\}\\
		& K:=\{k| N_k\in \mathcal{K}_G(N_j)\}=\{k_{1}, ... k_{n_{k}}\}
	\end{align*}
	for any $l=1,...,n$.\\
	The set $C$ contains all the indexes of the $n_c$ children of $N_n$ while the set $K$ contains all the $n_k$ indexes of the kins of $N_j$. The set $P_l$ contains all the $n_{p_l}$ parents (but $N_j$) of a generic node $N_l$.\\
	For any $i\neq j$, define
	\begin{align}\label{eq:def yc}
		\eps_i:=x_i-\sum_{k\neq j}H_{ij}(z)H_{jk}(z)x_k
			-\sum_{k\neq j}H_{ik}(z)x_k.
	\end{align}
	Note that, for any $i\neq j$,
	\begin{align}\label{eq:eps_i as in smart lemma}
		\eps_{i}=e_{i}+H_{ij}(z)e_{j}.
	\end{align}
	Also note that
	\begin{align}\label{eq: ej function of kins}
		e_j:=x_j-\sum_{k}H_{jk}(z)x_k,
	\end{align}
	The solution to the minimization problem
	\begin{align}
		&\arg\min_{q\in \text{tf-span}\{\eps_i\}_{i\neq j}}\|e_j-q\|^2=\\
		&\qquad =\arg\min_{\stackrel{q\in \text{tf-span}\{\eps_i\}_{i\neq j,0}}{\hat C_{jj}(z)\in~{}^{0}\mathcal{F}} }\|e_j-q\|^2+\|C_{jj}(z)e_0\|^2=\\
		&\qquad \arg\min_{q\in \text{tf-span}\{\eps_i\}_{i\neq j,0}}\|e_j-q\|
	\end{align}
	is given by
	\begin{align}
		\hat e_j = \sum_{\stackrel{i=1}{i\neq j}}^{n} \hat C_{ji}(z)\eps_{i},
	\end{align}
	where the process $e_0$ has played the role of a placeholder forcing $\hat C_{jj}=0$, without affecting the values of $\hat C_{ji}(z)$, for $i\neq j$.
	Now, let us consider the problem
	\begin{align}
		\arg\min_{q\in \text{tf-span}\{x_i\}_{i\neq j}}\|x_j-q\|.
	\end{align}
	From (\ref{eq: ej function of kins}), its solution $\hat x_j$ satisfies
	\begin{align}
		\hat x_j 	&= \sum_{k}H_{jk}(z)x_k+
			\arg\min_{q\in \text{tf-span}\{x_i\}_{i\neq j}}\|e_j-q\|=\\
			&= \sum_{k}H_{jk}(z)x_k+
			\arg\min_{q\in \text{tf-span}\{\eps_i\}_{i\neq j}}\|e_j-q\|
	\end{align}
	where the last equality has been obtained by using Lemma (\ref{lem:smart span}).
	Thus we have
	\begin{align}
		\hat x_j & = \sum_{i\neq j} W_{ji}x_i
			= \sum_{k}H_{jk}(z)x_k + \sum_{i\neq j} \hat C_{ji}\eps_{i}
	\end{align}
	Substituting the espression of $\eps_i$, $i\neq j$, as a function of $x_i$, $i \neq j$, the assertion is proven.
\end{IEEEproof}

Theorem \ref{thm:kinship explicit} gives the expression of the entry $W_{ji}$ of the topological filter as the sum of the three components $\hat C_{ji}$, $\hat P_{ji}$ and $\hat K_{ji}$.
The three components have a graphical interpretation:
$\hat C_{ji}(z)$ is a the contribution corresponding to the fact that $N_i$ is a child of $N_j$; $\hat P_{ji}(z)$ is the contribution present when $N_i$ is a parent of $N_j$; and $\hat K_{ji}(z)$ is a term present when there is a ``co-parent'' relation.
\begin{lem}\label{lem: three components}
	Given a well posed LDG, consider Equation (\ref{eq: child filter}), Equation (\ref{eq: parent filter}) and Equation (\ref{eq: kin filter}). Then
	\begin{itemize}
		\item	$\hat C_{ji}(z)\neq 0$ if and only if $N_i \in \mathcal{C}_G(N_j)$
		\item	$\hat P_{ji}(z)\neq 0$ if and only if $N_i \in \mathcal{P}_G(N_j)$
		\item	$\hat K_{ji}(z)\neq 0$ implies $N_i \in \mathcal{P}_G(\mathcal{C}_G(\mathcal{N}_j))\backslash N_j$.
	\end{itemize}
\end{lem}
\begin{IEEEproof}
	Evaluate
	\begin{align*}
		&\hat C_{j*}=\Phi_{e_j}H^{*}_{*j}(z)
			[\Phi_{e_j}(z)H_{*j}(z)H^{*}_{*j}(z)+\Phi_{e_{I_j}}(z)]^{-1}=\\
		&~=\Phi_{e_j}(z) H_{*j}^{*}(z)\Phi_{e_{I_j}}^{-1}(z)
			\left[
				\I+\Phi_{e_j}(z)H_{*j}(z)H_{*j}^*(z)\Phi_{e_{I_j}}(z)^{-1}
			\right]^{-1}=\\
		&~=\Phi_{e_j}(z) H_{*j}^{*}(z)\Phi_{e_{I_j}}^{-1}(z)
			\sum_{m=0}^{+\infty}
				\frac{(-1)^m}{m!}\|H_{*j}^*(z)\Phi_{e_{I_j}}^{-1}H_{*j}(z)\|^m\Phi_{e_j}(z)^m=\\
		&~=\Phi_{e_j}(z) H_{*j}^{*}(z)\Phi_{e_{I_j}}^{-1}(z)
		\frac{1}{1+\|H_{*j}^*(z)\Phi_{e_{I_j}}^{-1}(z)H_{*j}(z)\|\Phi_{e_j}(z)}.
	\end{align*}
	Thus $\hat C_{ji}\neq 0$ for any $N_i \in \mathcal{C}_G(N_j)$.
	Furthermore, it follows
	\begin{align}
		\|\hat C_{j*}(\w) H_{*j}(\w)\|<1,
	\end{align}
	implying, from (\ref{eq: parent filter}) , that
	\begin{align}
		\hat P_{ji}(\w)\neq 0.
	\end{align}
	for any $N_i \in \mathcal{P}_G(N_j)$.
\end{IEEEproof}
From Lemma \ref{lem: three components} Theorem \ref{thm:sparse wiener} can be recovered, since $\hat C_{ji}(z)=\hat P_{ji}(z)=\hat K_{ji}(z)=0$ for all the nodes $N_i$ that are not kin of $N_j$. However the only conditions according to which $W_{ji}=W_{ij}=0$ even though $\{N_i,N_j\} \in kin(G)$, making the link not detectable, are found.
\begin{theo}\label{thm:kinship null entry}
	Consider a well posed LDG $\mathcal{G}=(G,\mathcal{H},\mathcal{E})$ with $G=(V,A)$, $V=\{N_1,...,N_n\}$. $\{N_i, N_j\}\in kin(G)$ and $W_{ji}=W_{ij}=0$ if and only if one of the following condition is satisfied
	\begin{itemize}
		\item	$\hat C_{ji}(z)=-\hat K_{ji}(z)$ and $\hat P_{ij}(z)=-\hat K_{ij}(z)$
		\item   $\hat P_{ji}(z)=-\hat K_{ji}(z)$ and $\hat C_{ij}(z)=-\hat K_{ij}(z)$
	\end{itemize}

\end{theo}
\begin{IEEEproof}
	Indeed, since in a LDG $\mathcal{C}_G(N_j) \cap \mathcal{P}_G(N_j)=\emptyset$, it follows that, for any couple $(i,j)$, $\hat C_{ji}(z)=0$ or $\hat P_{ji}(z)=0$.
	The conditions, then, follow in a straightforward manner.
\end{IEEEproof}
Theorem \ref{thm:kinship null entry} shows that the condition under which a link present in the kin topology leads to two null entries ($W_{ji}=W_{ij}=0$) in the topological filter can be considered a pathological case. Indeed two simultaneous cancellations between the kin term and the other term (parent or child) must happen.
Moreover, if a node is simple, such a cancellation is not possible.


\section{Sparsity of causal filtering operators}\label{sec:causal}
First, we need to introduce the following lemma.
\begin{lem}\label{lem:zero span from generators}
	Let $x$ and $x_1,...,x_n$ be processes in the space $^0\mathcal{F}\Eps$.
	Define $I:=(1,....,n)$. Assume that $\Phi_{xx_I}(\w)=0$ for all $\w\in [-\pi, \pi]$.
	Then $<x,y>=0$ for all $y\in \text{tf-span}\{x_i\}_{i=1,...,n}$.
\end{lem}
\begin{IEEEproof}
	As $y\in \text{tf-span}\{x_i\}_{i=1,...,n}$, it follows that there exist $\alpha_i(z)\in {}^0\mathcal{F}$ such that
	\begin{align*}
		y=\sum_{i=1}^{n}\alpha_i(z)x_i=\alpha_I(z)x_I,
	\end{align*}
	where $\alpha_{I}(z)=[\alpha_{1}(z), ... , \alpha_{n}(z)]$ is a row vector of real-rational transfer functions.
	Then it follows that
	\begin{align*}
		<x,y>=\int_{-\pi}^{\pi} \Phi_{x x_I}(\w)\alpha_I(\w)^{*}=0.
	\end{align*}
\end{IEEEproof}
Now, a specific formulation of the causal Wiener filter is introduced for the defined spaces.
\begin{prop}\label{eq: my causal wiener}
	Let $x$ and $x_1,...,x_n$ be processes in the space $^0\mathcal{F}\Eps$.
	Define $I:=(1,....,n)$ and $X:=\text{c-tf-span}\{x_1,...,x_n\}$. Consider the problem
	\begin{align}\label{eq: cost general wiener (causal)}
		\inf_{q \in X} \|x-q\|^2.
	\end{align}
	Let $S(z)$ be the spectral factorization of $\Phi_{x_I}(\w)=S(\w)S^{*}(\w)$.
	If $\Phi_{x_I}(\w)>0$, for $\w\in[-\pi,\pi]$, the solution $\hat x^{(c)}\in X$ exists, is unique and has the form
	\begin{align*}
		\hat x^{(c)}=W^{(c)}(z)x_I
	\end{align*}
	where
	\begin{align*}
		W^{(c)}(z)=\{\Phi_{x x_{I}}(z)\Phi_{x_{I}}(z)^{-1}S(z)\}_{C}S^{-1}(z).
	\end{align*}
	Moreover $\hat x^{(c)}$ is the only element in $X$ such that, for any $q\in X$, it holds that
	\begin{align}\label{eq:perp hilbert projection thm (causal)}
		<x-\hat x^{(c)},q>=0.
	\end{align}
\end{prop}
\begin{IEEEproof}
	Observe that $W^{(c)}(z)$ is rational and causal. Let us prove that, for any causal $H(z)$, we have
	\begin{align}
		<x-\hat x^{(c)},H(z) x_{I}>=0
	\end{align}
	Define $r_{I}:=S^{-1}(z)x_{I}$, and observe that
	$S(z), S(z)^{-1}(z)\in~{}^0\mathcal{F}^{+}$,
	since $S(z)$ is the spectral factor of $\Phi_{x_I}(z)$ that is real-rational and has full rank on the unit circle. The signal $r_I$ is white since
	\begin{align}
		\Phi_{r_I}(z)=\mathcal{I}_{n}.
	\end{align}
	Define
	\begin{align}
		\hat x^{(c)}& =\{\Phi_{x r_I}(z)\Phi^{-1}_{r_I}(z)\}_{C}r_I\\
				& =\{\Phi_{x x_I}(z)\Phi^{-1}_{x_I}(z)S(z)\}_{C}S^{-1}(z)x_I.
	\end{align}
	observing that $\hat x^{(c)}\in X$ since it is obtained by composing transformations in $~{}^0\mathcal{F}^{+}$.
	Also observe that $\Phi_{x r_I}(z)\Phi^{-1}_{r_I}(z)$ is the non-causal Wiener filter estimating $x$ from $r_I$. Let the strictly anticausal component of such a filter part be
	\begin{align}
		\{\Phi_{x r_I}(z)\Phi^{-1}_{r_I}(z)\}_{SA}:=
			\Phi_{x r_I}(z)\Phi^{-1}_{r_I}(z)-\{\Phi_{x r_I}(z)\Phi^{-1}_{r_I}(z)\}_{C}.
	\end{align}
	Now, for a $H(z)\in~{}^0\mathcal{F}$, compute
	\begin{align}
		& <x-\hat x, Hx_I>=<x-W^{(c)}Sr_I x, HSr_I>=\\
		& <x-\Phi_{x r_I}(z)\Phi^{-1}_{r_I}(z)r_I+\{\Phi_{x r_I}(z)\Phi^{-1}_{r_I}(z)\}_{AC}r_I, HSr_I>=\\
		& <\{\Phi_{x r_I}(z)\Phi^{-1}_{r_I}(z)\}_{SA}r_I, HSr_I>=0.
	\end{align}
	where Lemma \ref{lem:zero span from generators} has been used to prove that
	$<x-\Phi_{x r_I}(z)\Phi^{-1}_{r_I}(z)r_I,H(z)S(z) r_I>=0$ and the last equality follows from the causality of $H(z)S(z)$, the strictly anti-causality of $\{\Phi_{x r_I}(z)\Phi^{-1}_{r_I}(z)\}_{SA}$ and the whiteness of $r_{I}$.
	From the Hilbert projection theorem, $\hat x^{(c)}$ is the unique process minimizing the cost (\ref{eq: cost general wiener (causal)}).
\end{IEEEproof}

The following theorem proves the sparsity of the causal Wiener filter stating that the causal Wiener filter estimating $x_{j}$ from the signals $x_{i}$, $i\neq j$, has non-zero entries corresponding to the kin signals of $x_{j}$.
\begin{theo}\label{thm:sparse causal wiener}
	Consider a well-posed, causal and topologically detectable LDG.
	Let $x_1,...,x_n\in {}^0\mathcal{F}\Eps$ be the signals associated with the $n$ nodes of its graph.
	Define $X_j=\text{c-tf-span}\{x_i\}_{i\neq j}$.
	Consider the problem of approximating the signal $x_j$ with an element $\hat x_j\in X_j$, as defined below
	\begin{align}
		\min_{\hat x_j\in X_j}
		\left\|x_j-\hat x_j \right\|^2.
	\end{align}
	Then the optimal solution $\hat x_j$ exists, is unique and
	\begin{align}
		\hat x_j =\sum_{i\neq j} W_{ji}(z)x_i
	\end{align}
	where $W_{ji}(z)\neq 0$ implies $(N_i, N_j) \in kin(G)$.
\end{theo}
\begin{IEEEproof}
	For any $i\neq j$, define $\eps_i$ as in (\ref{eq:def yc}) and observe that $\eps_i$ can be represented as in (\ref{eq:eps_i as in smart lemma}).
	Also note that
	\begin{align}\label{eq: e_j causal}
		e_j:=x_j-\sum_{i}H_{ji}(z)x_i.
	\end{align}
	Consider $\hat e_j$ defined as
	\begin{align}
		\hat e_j
		:=\arg\min_{q\in \text{c-tf-span}\{\eps_i\}_{i\neq j}}\|e_j-q\|
		=\sum_{i\neq j} C^{(c)}_{ji}(z)\eps_{i}
	\end{align}
	where the transfer functions $C^{(c)}_{ji}(z)$ are given by the causal Wiener filter estimating $e_j$ from $\{\eps_{i}\}_{i\neq j}$. Notice that, by (\ref{eq:eps_i as in smart lemma}), $C^{(c)}_{ji}(z)$ is equal to zero if $N_i$ is not a child of $N_j$.
	Now, let us consider the optimization problem
	\begin{align}
		\hat x_{j}:=\arg\min_{q\in \text{c-tf-span}\{x_i\}_{i\neq j}}\|x_j-q\|=\sum_{i\neq j}W_{ji}(z)x_i
	\end{align}
	where $W_{ji}(z)$ are the entries of the causal Wiener filter solving it.
	Its solution $\hat x_j$ satisfies also
	\begin{align}
		\hat x_j 	&= \sum_{i\neq j}H_{ji}(z)x_i+
			\arg\min_{q\in \text{c-tf-span}\{x_i\}_{i\neq j}}\|e_j-q\|=\\
			&= \sum_{i}H_{ji}(z)x_i+
			\arg\min_{q\in \text{c-tf-span}\{\eps_i\}_{i\neq j}}\|e_j-q\|
	\end{align}
	where the first equality derives from (\ref{eq: e_j causal}) and the last one has been obtained by using Lemma \ref{lem:smart span}.
	Thus we have
	\begin{align}
		\hat x_j & = \sum_{i\neq j} W_{ji}x_i
			= \sum_{i}H_{ji}(z)x_i + \sum_{i\neq j} C_{ji}\eps_{i}
	\end{align}
	Substituting the espression of $\eps_i$, $i\neq j$, as a function of $x_i$, $i \neq j$, the assertion is proven.
\end{IEEEproof}

The following theorem proves the sparsity of the one step prediction operator (or granger-causal operator).
Granger-causality is a wide-spread technique in econometrics to test the causal dependence of time series.
If the stronger hypothesis of strictly causal transfer functions $H_{ji}(z)$ is met Granger-causality provides an exact reconstruction of parent-child links in a LDG.
\begin{theo}\label{thm:sparse granger}
	Consider a well-posed, strictly causal and topologically detectable LDG.
	Let $x_1,...,x_n\in {}^0\mathcal{F}\Eps$ be the signals associated with the $n$ nodes of its graph.
	Define $X_j=\text{c-tf-span}\{x_i\}$.
	Consider the problem of approximating the signal $z~x_j$ with an element $\hat x_j\in X_j$, as defined below
	\begin{align}
		\min_{\hat x_j\in X_j}
		\left\|z x_j-\hat x_j \right\|^2.
	\end{align}
	Then the optimal solution $\hat x_j$ exists, is unique and
	\begin{align}
		\hat x_j =\sum_{i=1}^{n} W_{ji}(z)x_i
	\end{align}
	where $W_{ji}(z)\neq 0$ implies $i=j$ or $N_i$ is a parent of $N_i$.
\end{theo}
\begin{IEEEproof}
	For any $i\neq j$, define $\eps_i$ as in (\ref{eq:def yc}) and observe that $\eps_i$ can be represented as in (\ref{eq:eps_i as in smart lemma}). Also define $\eps_j:=e_j$.
	Note that
	\begin{align}\label{eq: e_j granger}
		e_j:=x_j-\sum_{i}H_{ji}(z)x_i.
	\end{align}
	Consider the minimization problem
	\begin{align}
		\hat e_{j}:=\arg\min_{q\in \text{c-tf-span}\{\eps_i\}_{i}}\|ze_j-q\|= \sum_{i\neq j} C^{(g)}_{ji}(z)\eps_{i}
	\end{align}
	where the transfer functions $C^{(g)}_{ji}(z)$ are elements of ${}^{0}\mathcal{F}$.
	We have that $C^{(g)}_{ji}(z)=0$ for any $i\neq j$. Indeed, it holds that
	\begin{align}
		&\arg\min_{q\in \text{c-tf-span}\{\eps_i\}_{i=1}^{n}}\|ze_j-q\|=
		\arg\min_{q\in \text{c-tf-span}\{e_i\}_{i=1}^{n}}\|ze_j-q\|=\\
		&\qquad = \arg\min_{q\in \text{c-tf-span}\{e_j\}}\|ze_j-q\|
	\end{align}
	and, since $\Phi_{e_ie_j}(\w)=0$ for $i\neq j$.
	Conversely, from (\ref{eq: my causal wiener}), we find$C^{(g)}_{jj}(z)=\{zS_{j}(z)\}_{C}z^{-1}S_j(z)$ where $S_j(z)$ is the spectral factor of $e_j$.
	Now, let us consider the problem
	\begin{align}
		\arg\min_{q\in \text{c-tf-span}\{x_i\}_{i\neq j}}\|zx_j-q\|.
	\end{align}
	Its solution $\hat x_j$ is
	\begin{align*}
		\hat x_j 	&= \sum_{k}zH_{jk}(z)x_k+
			\arg\min_{q\in \text{c-tf-span}\{x_i\}_{i}}\|ze_j-q\|=\\
			&= \sum_{k}zH_{jk}(z)x_k+C^{(g)}_{jj}(z)e_j\\
			&= C^{(g)}_{jj}(z)x_j+\sum_{k\neq j} [zH_{jk}(z)-C^{(g)}_{jj}(z)H_{jk}(z)]x_k.
	\end{align*}
	This proves the assertion.
\end{IEEEproof}

\section{A reconstruction algorithm}\label{sec:algorithm}
The previous section provides theoretical results allowing for the reconstruction of a topology via Wiener filtering. It needs to be stressed that even in the case of sparse graphs, the reconstruction of the kinship topology can be considered a practical solution. The reasons are two-fold.
In many situations it is possible to measure the outputs of many nodes, while it is important to identify a reduced number of possible interconnections among those nodes. For example, DNA-microarrays are devices that allow the measurement of gene expression of a cell of a cell.
Such data can be useful in understanding which genes interact together and realize a specific metabolic pathway and how they are related \cite{RavSom02}.
Indeed, a cell can express tens of thousands of genes while only a few tens must be involved in a gene regulatory network. The possibility of reducing the set of involved genes to test the presence of actual interactions with targeted experiments is of significant importance \cite{MarPel08b}.
Analogously, in Finance, quantifying the strongest interconnections among a set of market stocks can suggest good strategies to balance a given portfolio \cite{ManSta00}, thus it is important to have a quantitative tool to group together different stocks or, at least, to detect a limited set of possible dynamical connections.
Similar problems are also faced in neuroscience in order to understand neural interconnections \cite{BroDin04}.
A second reason why the presented analysis is of practical importance is that as a byproduct of the reconstruction an optimal model for the node dynamics which can be used for smoothing procedures (in the case non-causal Wiener filter are derived) or predictive ones (in the strictly causal case which is out of the scope of this article) can be obtained.

The following algorithm is the natural pseudocode implementation of the reconstruction technique that was developed in the previous section.
\subsection*{\tt Reconstruction algorithm}
\begin{itemize}
	\item[{\tt 0.}] Initialize the set of edges $A=\{ \}$
	\item[{\tt 1.}] For any signal $x_j$
	\item[{\tt 2.}]\quad solve the minimization problem (\ref{eq:min problem})
	\item[{\tt 3.}]\quad For any $W_{ji}(z) \ncong 0$
	\item[{\tt 4.}]\quad\quad add $\{N_i,N_j\}$ to $A$
	\item[{\tt 5.}]\quad end
	\item[{\tt 7.}] end
	\item[{\tt 8.}] return $A$
\end{itemize}

Step {\tt 2} can be computed efficiently for a large number of signals using techniques based on Gram-Schmidt orthogonalization as those described in \cite{GolRee98}.
When checking the condition at step {\tt 3}, the condition $W_{ji} \ncong 0$ can be implemented as $\|W_{ji}(z)\|>\sigma_{thr}(x_i,x_j)$ where $\sigma_{thr}(x_i,x_j)$ is a threshold that might depend on both the signals $x_i$ and $x_j$. Moreover, a proper norm must be defined on the space of transfer functions.\\

\section{Discussion on the robustness of the reconstruction}\label{sec:corrupted}
In the previous sections an ideal scenario was analyzed, where the node signals are  affected by the process noises $e_{j}$ which are actively driving the network.
This leads to the result that the Wiener filter is a useful tool for identifying the network topology providing the topology of the smallest self-kin network that contains the original system.
In most pratical cases the signals are not accessible but are also affected by measurement noise.
Unfortunately, the presence of measurement noise destroys the sparsity property of the Wiener filter. 
However, it will be shown that when the measurement noise is known to be sufficiently small, the introduction of a threshold in the algorithm again guarantees the reconstruction of the topology of the minimal self-kin network.
In different terms, the Wiener filter property presented in order to reconstruct a topology offers a certain degree of robustness that will be quantified.\\
\begin{define}
	A Corrupted LDG (CLDG) is a pair $(\mathcal{G},\eta_I)$ where $\mathcal{G}=(G,\mathcal{H},e_{I})$ is a LDG with dynamics $X=\{x_1 ,..., x_{n}\}$ and $\eta_I=\{\eta_1,...,\eta_n\}$ is a set of noises with the property that they are mutually not correlated and not correlated with the signals $\{e_{j}\}\in \mathcal{E}$, either.
	The output of the CLDG as the set $Y$ of signals is defined as
	\begin{align*}
		y_j=x_j+\eta_j.
	\end{align*}
	Introducing an auxiliary signal $y_0$ which is not correlated with any $y_j$ $(j=1,...,n)$,
	the additively corrupted topological filter of $\mathcal{G}$ with respect to the measures $y_{j}$ is defined as
	\begin{align}
		\hat W(z)=\left(\begin{array}{c}
				\Phi_{y_1 y_{I_{1}}}(z)\Phi_{y_{I_1} y_{I_1}}^{-1}(z)\\
				\vdots \\
				\Phi_{y_n x_{I_{n}}}(z)\Phi_{y_{I_n} y_{I_n}}^{-1}(z)
			\end{array}\right).
	\end{align}
\end{define}
The following important lemma about the inverse of a matrix is recalled.
\begin{lem}\label{lem:inverse lemma}
	If $Q$ and $Q+\Delta$ are two invertible matrices the following equality holds
	\begin{align}
		(Q+\Delta)^{-1}-Q^{-1}=-Q^{-1}(Q^{-1}+\Delta^{-1})^{-1}Q^{-1}.
	\end{align}
\end{lem}
The following result establishes a bound on the difference between the topological filter $W(z)$ which can be obtained by measuring directly the dynamics of a LDG and the corrupted topological filter $\hat W(z)$ obtained by measuring its output.
\begin{theo}\label{thm:bound on corrupted}
	Consider a CLDG $(\mathcal{G},\eta)$. Let $W(z)$ be the topological filter of $\mathcal{G}$ and let $\hat W(z)$ be its corrupted topological filter.
	If
	\begin{align}
		\frac{1}{\|\Phi_{y_{I_j} y_{I_j}}^{-1}(z)\|}-\|\Phi_{\eta}(z)\|>0,
	\end{align}
	then, for any $z\in \mathbf{C}$,
	\begin{align}
		\|\hat W_{j}(z) - W_{j}(z) \|&\leq 
		\frac{\|\Phi_{y_{j} y_{I_j}}(z) \|\|\Phi_{y_{I_j} y_{I_j}}(z)^{-1} \| \|\Phi_{\eta}(z) \|}
		{\left( \frac{1}{\|\Phi_{y_{I_j} y_{I_j}}^{-1}(z)\|}-\|\Phi_{\eta}(z)\|
		\right)}.
	\end{align}
\end{theo}
\begin{IEEEproof}
	First note that
	\begin{align}
		\Phi_{y_{j} y_{I_j}}(z)=\Phi_{x_{j} x_{I_j}}(z),
	\end{align}
	thus
	\begin{align}
		\hat W_{j}(z) - W_{j}(z)=\Phi_{y_{j} y_{I_j}}(z)
			\left[
			\Phi_{y_{I_j} y_{I_j}}^{-1}(z)-\Phi_{x_{I_j} x_{I_j}}^{-1}(z)
			\right]=\\
			\Phi_{y_{j} y_{I_j}}(z)
			\left[
			\Phi_{y_{I_j} y_{I_j}}^{-1}(z)-
				\left(
					\Phi_{y_{I_j} y_{I_j}}(z)-\Phi_{\eta}(z)
				\right)^{-1}
			\right].
	\end{align}
	By applying Lemma~\ref{lem:inverse lemma}, it follows that
	\begin{align}
		&Q^{-1}-(Q-\Delta)^{-1}=Q^{-1}(Q^{-1}-\Delta^{-1})^{-1}Q^{-1}=\\
		& =Q^{-1}(\Delta-Q)^{-1}\Delta,
	\end{align}
	which implies
	\begin{align}
		&\|Q^{-1}-(Q-\Delta)^{-1}\|\leq \|Q^{-1}\|\|\Delta\|\|(Q-\Delta)^{-1}\|=\\
		&=\frac{\|Q^{-1}\|\|\Delta\|}{\min sp(Q-\Delta)}
			\leq\frac{\|Q^{-1}\|\|\Delta\|}{\min sp(Q)-\max sp(\Delta)}=\\
		&=\frac{\|Q^{-1}\|\|\Delta\|}{\frac{1}{\|Q^{-1}\|}-\|\Delta\|}
	\end{align}

\end{IEEEproof}
An immediate way to apply Theorem~\ref{thm:bound on corrupted} is when a function $D(z)$ bounding $\|\Phi_{\eta}(z)\|$ is known as illustrated by the following corollary.
\begin{cor}
	Consider a CLDG $(\mathcal{G},\eta)$. Let $\hat W(z)$ be its corrupted topological filter.
	Assume that there exists a real valued function $D(z)$ such that
	\begin{align}
		\|\Phi_{\eta}(z)\|<D(z).
	\end{align}
	If, for some $z\in \mathbf{C}$,
	\begin{align}
		\frac{1}{\|\Phi_{y_{I_j} y_{I_j}}^{-1}(z)\|}-\|D(z)\|>0
	\end{align}
	and
	\begin{align}
		\|\hat W_{ji}(z)\|\geq 
		\frac{\|\Phi_{y_{j} y_{I_j}}(z) \| \|\Phi_{\eta}(z) \|}
		{\left( \frac{1}{\|\Phi_{y_{I_j} y_{I_j}}^{-1}(z)\|}-\|\Phi_{\eta}(z)\|
		\right)},
	\end{align}
	then
	\begin{align}
		(N_j, N_i)\in kin(G).
	\end{align}
\end{cor}
\begin{IEEEproof}
	By contradiction, assume that $(N_j, N_i)\notin kin(G)$. Then $W_{ji}(z)=0$.
	By Theorem~\ref{thm:bound on corrupted} we have 
	\begin{align}
		\|\hat W_{ji}(z)\|
			&\leq 
		\frac{\|\Phi_{y_{j} y_{I_j}}(z) \| \|\Phi_{\eta}(z) \|}
		{\left( \frac{1}{\|\Phi_{y_{I_j} y_{I_j}}^{-1}(z)\|}-\|\Phi_{\eta}(z)\|
		\right)}\\
			& <
		\frac{\|\Phi_{y_{j} y_{I_j}}(z) \| \|\Phi_{\eta}(z) \|}
		{\left( \frac{1}{\|\Phi_{y_{I_j} y_{I_j}}^{-1}(z)\|}-\|D(z)\|
		\right)}
	\end{align}
	which is a contradiction.
\end{IEEEproof}

Note that the inequality of Theorem~\ref{thm:bound on corrupted} needs to hold for any $z$, so, in this sense, it constitutes a sharp criterion to detect the presence of a link in a CLDG.

\section{Numerical Examples}\label{sec:example}
In this section illustrative applications of the theoretical results are provided.
It is worth observing that results were developed for the general class of linear models.
Indeed, no assumptions were made on the order and causality property of considered transfer functions.

\subsection{Self-kin network}\label{sec:self-kin example}
\begin{figure}[bth]
	\centering
	\begin{tabular}{cc}
		\includegraphics[width=0.45\columnwidth]{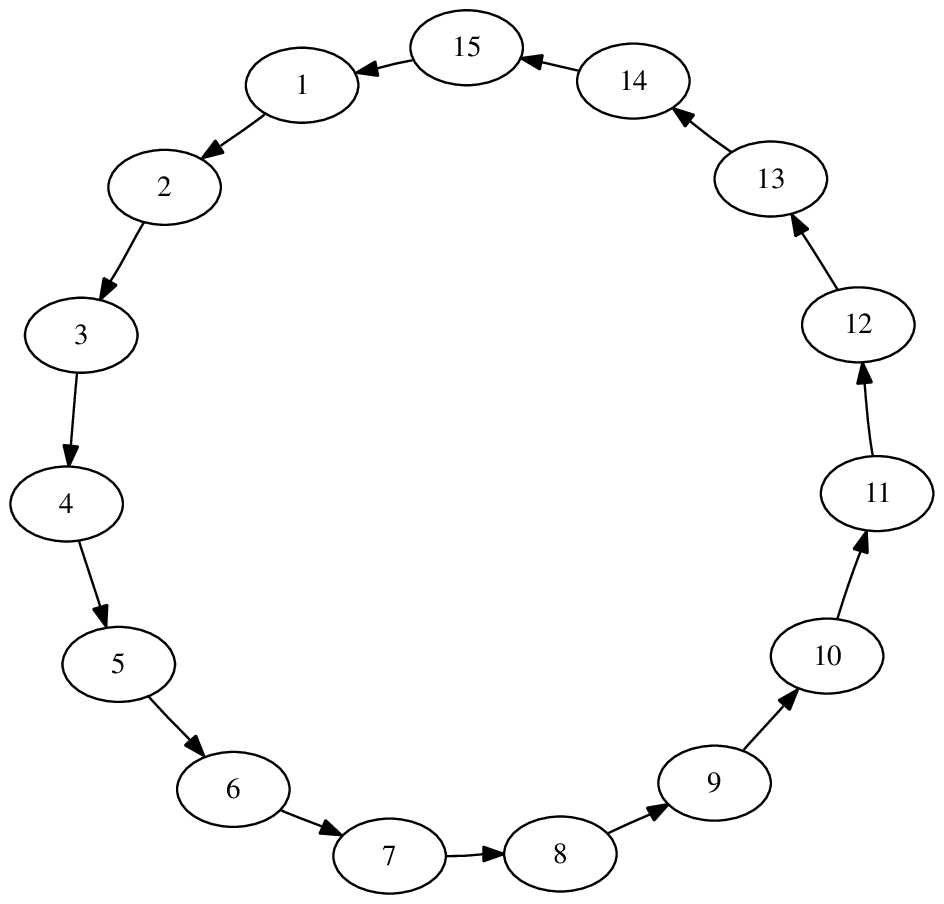} &
		\includegraphics[width=0.45\columnwidth]{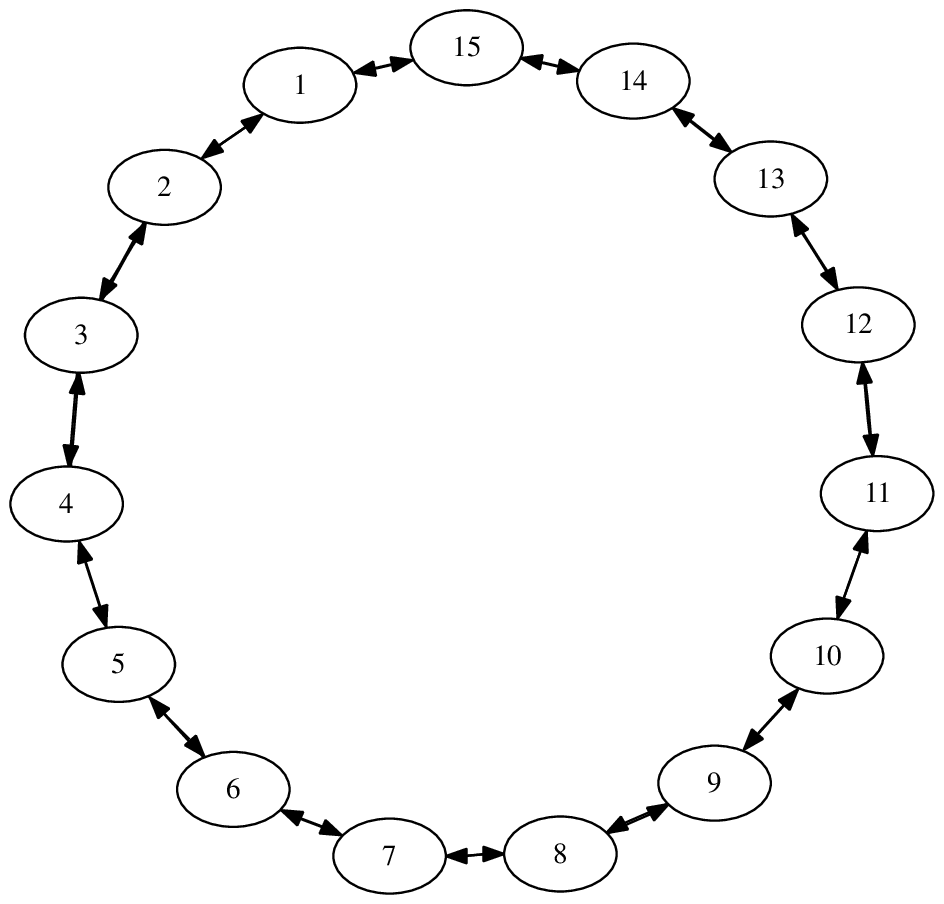} \\
		(a) & (b)
	\end{tabular}
	\caption{A ring network of $15$ nodes as the one considered in Example ~\ref{sec:self-kin example} (a).\label{fig:true_ring}
	The reconstructed topology of Example ~\ref{sec:self-kin example}. Every single links has been detected and, since the network is self-kin, the topology does not contain any spurious link (b). \label{fig:rec_ring}
	}
\end{figure}

Consider a ring network of $15$ nodes where the dynamics of the links is given by $5$-th order FIR filters and the noise power is the same on every node. The network topology is provided in Figure~\ref{fig:true_ring}
The network was simulated for $1000$ steps and an implementation of the developed algorithm was applied to the data providing the topology of Figure~\ref{fig:rec_ring}(b).

\subsection{Generic network}\label{sec:generic example}
\begin{figure}[bth]
	\centering
	\begin{tabular}{cc}
		\includegraphics[width=0.45\columnwidth]{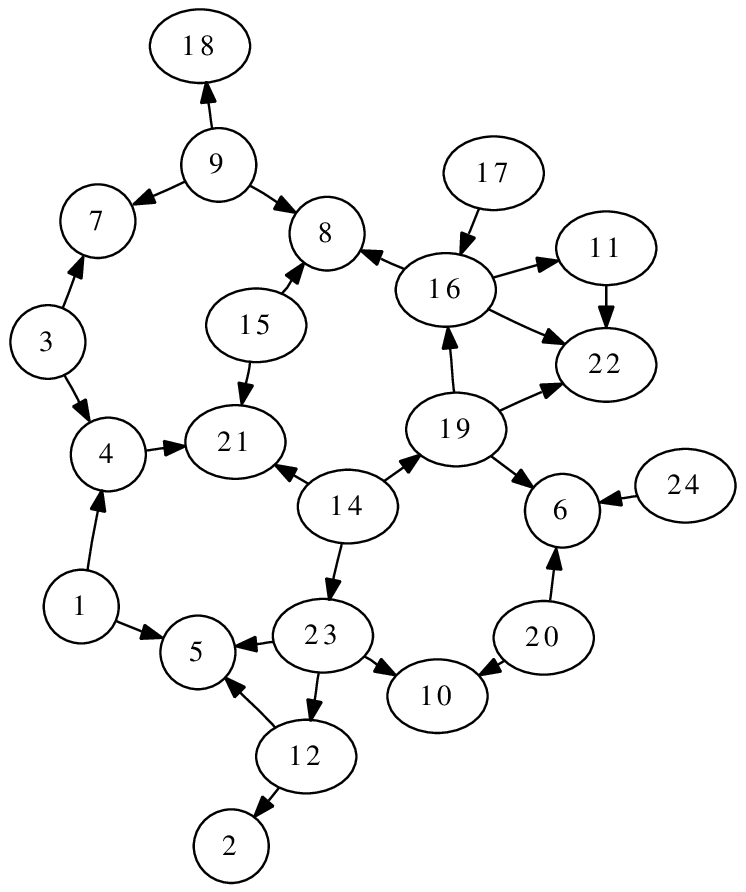} &
		\includegraphics[width=0.45\columnwidth]{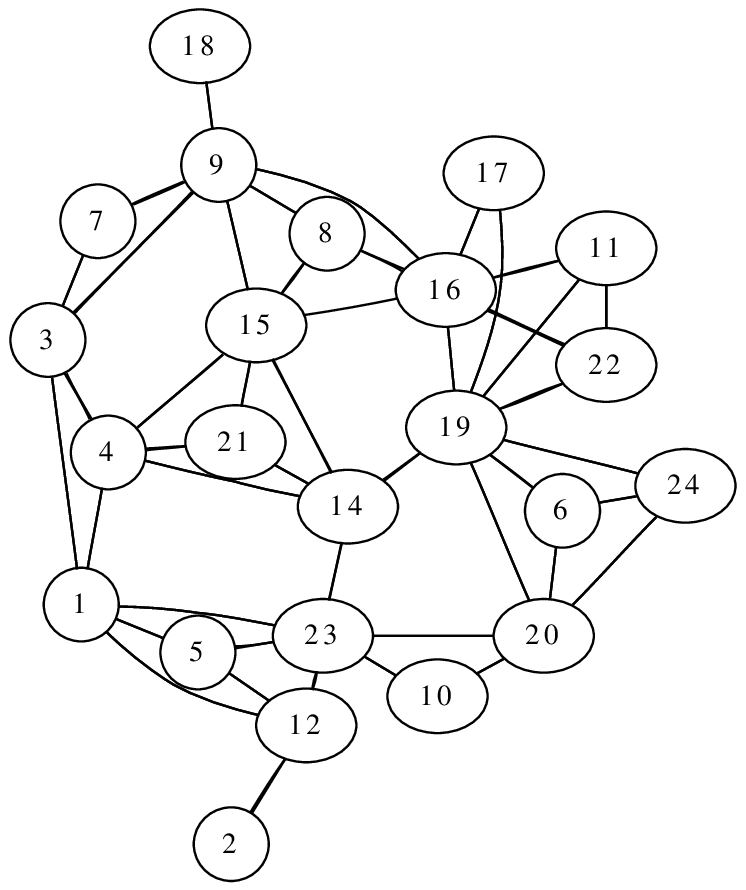}\\
		(a) & (b)
	\end{tabular}
	\caption{
		A network of $24$ nodes as the one considered in Example ~\ref{sec:generic example} (a) and the reconstructed topology (b).
		Every single links has been detected, but, since the network is not self-kin, as expected, the topology contains the additional links between the ``kins''.
		\label{fig:true_net}
		\label{fig:rec_net}
	}
\end{figure}

Consider a network of $24$ nodes for $1000$ steps as reported in Figure~\ref{fig:true_net}(a) and the reconstructed topology obtained is depicted in Figure~\ref{fig:rec_net}(b).

\section{Conclusions}
This work has illustrated a simple but effective procedure to identify 
the general structure of a network of linear dynamical systems.
To the best knowledge of the authors, a general analytical formulation of the problem of identifying a network is not tackled in scientific literature yet.
The approach followed is based on Wiener Filtering in order to detect the  existing links of a network.
When the topology of the original graph is described by a self-kin network, the method developed guarantees an exact reconstruction. Self-kin networks provide a non-trivial class of networks since they allow the presence of loops, nodes with multiple inputs and lack of connectivity.
Moreover, the paper also provides results about general networks. It is shown that, for a general graph, the developed procedure reconstructs the topology of the smallest self-kin graph containing the original one. Thus, the method is optimal in this sense.
Numerical examples illustrate the correctness and also the reliability of the identification technique.

\section*{Acknowlgements}
Donatello Materassi wants to thank Lipeng Ning for the useful discussions they had on the topics developed in this paper.

\bibliography{../../bibliography/topident,../../bibliography/l0problem,../../bibliography/control}

\begin{thebibliography}{10}
\providecommand{\url}[1]{#1}
\csname url@samestyle\endcsname
\providecommand{\newblock}{\relax}
\providecommand{\bibinfo}[2]{#2}
\providecommand{\BIBentrySTDinterwordspacing}{\spaceskip=0pt\relax}
\providecommand{\BIBentryALTinterwordstretchfactor}{4}
\providecommand{\BIBentryALTinterwordspacing}{\spaceskip=\fontdimen2\font plus
\BIBentryALTinterwordstretchfactor\fontdimen3\font minus
  \fontdimen4\font\relax}
\providecommand{\BIBforeignlanguage}[2]{{%
\expandafter\ifx\csname l@#1\endcsname\relax
\typeout{** WARNING: IEEEtran.bst: No hyphenation pattern has been}%
\typeout{** loaded for the language `#1'. Using the pattern for}%
\typeout{** the default language instead.}%
\else
\language=\csname l@#1\endcsname
\fi
#2}}
\providecommand{\BIBdecl}{\relax}
\BIBdecl

\bibitem{CziBar99}
A.~Czir\'{o}k, A.-L. Barab\'{a}si, and T.~Vicsek, ``Collective motion of
  self-propelled particles: Kinetic phase transition in one dimension,''
  \emph{Physical Review Letters}, vol.~82, no.~1, pp. 209+, January 1999.

\bibitem{LevRap00}
H.~Levine, W.-J. Rappel, and I.~Cohen, ``Self-organization in systems of
  self-propelled particles,'' \emph{Physical Review E}, vol.~63, no.~1, pp.
  017\,101+, December 2000.

\bibitem{FaxMur04}
A.~Fax and R.~M. Murray, ``Information flow and cooperative control of vehicle
  formations,'' \emph{IEEE Transactions on Automatic Control}, vol.~49, no.~9,
  pp. 1465--1476, 2004.

\bibitem{LiuYad08}
J.~Liu, V.~Yadav, H.~Sehgal, J.~M. Olson, H.~Liu, and N.~Elia, ``Phase
  transitions on fixed connected graphs and random graphs in the presence of
  noise,'' \emph{IEEE Transactions on Automatic Control}, vol.~53, p. 1817,
  2008.

\bibitem{Gee05}
D.~Geer, ``Industry trends: Chip makers turn to multicore processors,''
  \emph{Computer}, vol.~38, no.~5, pp. 11--13, 2005.

\bibitem{KadBhu08}
S.~C. Kadu, M.~Bhushan, and R.~Gudi, ``Optimal sensor network design for
  multirate systems,'' \emph{Journal of Process Control}, vol.~18, no.~6, pp.
  594 -- 609, 2008.

\bibitem{NayRos07}
M.~Naylora, L.~Roseb, and B.~Moyle, ``Topology of foreign exchange markets
  using hierarchical structure methods,'' \emph{Physica A}, vol. 382, pp.
  199--208, 2007.

\bibitem{ManSta00}
R.~Mantegna and H.~Stanley, \emph{An Introduction to Econophysics: Correlations
  and Complexity in Finance}.\hskip 1em plus 0.5em minus 0.4em\relax Cambridge
  UK: Cambridge University Press, 2000.

\bibitem{EisSpe98}
M.~Eisen, P.~Spellman, P.~Brown, and D.~Botstein, ``Cluster analysis and
  display of genome-wide expression patterns,'' \emph{Proc. Natl. Acad. Sci.
  USA}, vol.~95, no.~25, pp. 14\,863--8, 1998.

\bibitem{RavSom02}
E.~Ravasz, A.~Somera, D.~Mongru, Z.~Oltvai, and A.~Barabasi, ``Hierarchical
  organization of modularity in metabolic networks,'' \emph{Science}, vol. 297,
  p. 1551, 2002.

\bibitem{DelNin08}
D.~Del~Vecchio, A.~Ninfa, and E.~Sontag, ``Modular cell biology: Retroactivity
  and insulation,'' \emph{Nature Molecular Systems Biology}, vol.~4, p. 161,
  2008.

\bibitem{BroDin04}
A.~Brovelli, M.~Ding, A.~Ledberg, Y.~Chen, R.~Nakamura, and S.~L. Bressler,
  ``Beta oscillations in a large-scale sensorimotor cortical network:
  directional influences revealed by {G}ranger causality.'' \emph{Proc Natl
  Acad Sci USA}, vol. 101, no.~26, pp. 9849--9854, June 2004.

\bibitem{BunUrb00}
A.~Bunn, D.~Urban, and T.~Keitt, ``Landscape connectivity: A conservation
  application of graph theory,'' \emph{Journal of Environmental Management},
  vol.~59, no.~4, pp. 265 --278, 2000.

\bibitem{UrbKei01}
D.~Urban and T.~Keitt, ``Landscape connectivity: A graph-theoretic
  perspective,'' \emph{Ecology}, vol.~82, no.~5, pp. 1205--1218, 2001.

\bibitem{BaiMon06}
J.-S. Bailly, P.~Monestiez, and P.~Lagacherie, ``Modelling spatial variability
  along drainage networks with geostatistics,'' \emph{Mathematical Geology},
  vol.~38, no.~5, pp. 515--539, 2006.

\bibitem{MonBai05}
P.~Monestiez, J.-S. Bailly, P.~Lagacheriec, and M.~Voltz, ``Geostatistical
  modelling of spatial processes on directed trees: Application to fluvisol
  extent,'' \emph{Geoderma}, vol. 128, pp. 179--191, 2005.

\bibitem{InnMat08}
G.~Innocenti and D.~Materassi, ``A modeling approach to multivariate analysis
  and clusterization theory,'' \emph{Journal of Physics A}, vol.~41, no.~20, p.
  205101, 2008.

\bibitem{BocLat06}
S.~Boccaletti, V.~Latora, Y.~Moreno, M.~Chavez, and D.~U. Hwang, ``Complex
  networks: Structure and dynamics,'' \emph{Physics Reports}, vol. 424, no.
  4-5, pp. 175--308, February 2006.

\bibitem{GirNew02}
M.~Girvan and M.~E.~J. Newman, ``Community structure in social and biological
  networks,'' \emph{Proceedings of the National Academy of Sciences}, vol.~99,
  no.~12, 2002.

\bibitem{NewGir04}
M.~E.~J. Newman and M.~Girvan, ``Finding and evaluating community structure in
  networks,'' \emph{Physical Review E}, vol.~69, no.~2, 2004.

\bibitem{ZhaLiu07}
H.~Zhang, Z.~Liu, M.~Tang, and P.~Hui, ``An adaptative routing strategy for
  packet delivery in complex networks,'' \emph{Physics Letters A}, vol. 364,
  pp. 177--182, 2007.

\bibitem{Olf07}
R.~Olfati-Saber, ``Distributed kalman filtering for sensor networks,'' in
  \emph{Proc. of IEEE CDC}, New Orleans, 2007, pp. 5492--5498.

\bibitem{SchRib08}
I.~Schizas, A.~Ribeiro, and G.~Giannakis, ``Consensus in ad hoc {WSNs} with
  noisy links-part i: Distributed estimation of deterministic signals,''
  \emph{IEEE Trans. on Signal Processing}, vol.~56, no.~1, 2008.

\bibitem{Kol09}
E.~Kolaczyk, \emph{Statistical Analysis of Network Data: Methods and
  Models}.\hskip 1em plus 0.5em minus 0.4em\relax Berlin, Germany:
  Springer-Verlag, 2009.

\bibitem{Tim07}
M.~Timme, ``Revealing network connectivity from response dynamics,''
  \emph{Phys. Rev. Lett.}, vol.~98, no.~22, p. 224101, 2007.

\bibitem{BocIva07}
S.~{Boccaletti}, M.~{Ivanchenko}, V.~{Latora}, A.~{Pluchino}, and
  A.~{Rapisarda}, ``{Detecting complex network modularity by dynamical
  clustering},'' \emph{Phys. Rev. E}, vol.~75, p. 045102, 2007.

\bibitem{NapSau08}
D.~Napoletani and T.~Sauer, ``Reconstructing the topology of sparsely connected
  dynamical networks,'' \emph{Phys. Rev. E}, vol.~77, p. 026103, 2008.

\bibitem{OzeUzu08}
M.~Ozer and M.~Uzuntarla, ``Effects of the network structure and coupling
  strength on the noise-induced response delay of a neuronal network,''
  \emph{Physics Letters A}, vol. 375, pp. 4603--4609, 2008.

\bibitem{MicSok57}
C.~Michener and R.~Sokal, ``A quantitative approach to a problem of
  classification,'' \emph{Evolution}, vol.~11, pp. 490--499, 1957.

\bibitem{FreHar02}
R.~Freckleton, P.~Harvey, and M.~Pagel, ``Phylogenetic analysis and comparative
  data: A test and review of evidence,'' \emph{American Naturalist}, vol. 160,
  pp. 712--726, 2002.

\bibitem{Die06}
R.~Diestel, \emph{Graph Theory}.\hskip 1em plus 0.5em minus 0.4em\relax Berlin,
  Germany: Springer-Verlag, 2006.

\bibitem{LedWol04}
O.~Ledoit and M.~Wolf, ``Honey, i shrunk the covariance matrix,'' \emph{Journal
  of. Portfolio Management}, pp. 110--119, 2004.

\bibitem{TumLil08}
M.~Tumminello, F.~Lillo, and R.~N. Mantegna, ``Shrinkage and spectral filtering
  of correlation matrices: a comparison via the {K}ullback-{L}eibler
  distance,'' \emph{Acta Physica Polonica B}, pp. 4079--4088, 2008.

\bibitem{TotKer09}
B.~Toth and J.~Kertesz, ``Accurate estimator of correlations between
  asynchronous signals,'' \emph{Physica A}, pp. 1696--1705, 2009.

\bibitem{MatInn09}
D.~Materassi and G.~Innocenti, ``Unveiling the connectivity structure of
  financial networks via high-frequency analysis,'' \emph{Physica A:
  Statistical Mechanics and its Applications}, vol. 388, no.~18, pp.
  3866--3878, June 2009.

\bibitem{MatInn08}
------, ``Topological identification in networks of dynamical systems,'' in
  \emph{Proc. of IEEE CDC}, Cancun (Mexico), December 2008.

\bibitem{GetFri02}
L.~Getoor, N.~Friedman, B.~Taskar, and D.~Koller, ``Learning probabilistic
  models of relational structure,'' \emph{Journal of Machine Learning
  Research}, vol.~3, pp. 679--707, 2002.

\bibitem{FriKol03}
N.~Friedman and D.~Koller, ``Being bayesian about network structure: A bayesian
  approach tostructure discovery in bayesian networks,'' \emph{Machine
  Learning}, vol.~50, pp. 95--126, 2003.

\bibitem{Pea88}
J.~Pearl, \emph{{Probabilistic reasoning in intelligent systems: networks of
  plausible inference}}.\hskip 1em plus 0.5em minus 0.4em\relax Morgan
  Kaufmann, 1988.

\bibitem{MarPel08a}
D.~Marinazzo, M.~Pellicoro, and S.~Stramaglia, ``Kernel method for nonlinear
  {G}ranger causality,'' \emph{Phys. Rev. Lett.}, vol. 100, p. 144103, 2008.

\bibitem{Gra69}
C.~Granger, ``Investigating causal relations by econometric models and
  cross-spectral methods,'' \emph{Econometrica}, vol.~37, pp. 424--438, 1969.

\bibitem{CanWak08}
E.~Cand\`{e}s, M.~Wakin, and S.~Boyd, ``Enhancing sparsity by reweighted l1
  minimization,'' \emph{Journal of Fourier Analysis and Applications}, vol.~14,
  pp. 877--905, 2008.

\bibitem{CanTao05}
E.~J. Cand{\`e}s and T.~Tao, ``Decoding by linear programming,'' \emph{IEEE
  Transactions on Information Theory}, vol.~51, no.~12, pp. 4203--4215, 2005.

\bibitem{JovBam05}
M.~R. Jovanovi\'c and B.~Bamieh, ``On the ill-posedness of certain vehicular
  platoon control problems,'' \emph{IEEE Transactions on Automatic Control},
  vol.~50, no.~9, pp. 1307--1321, September 2005.

\bibitem{Lue69}
D.~G. Luenberger, \emph{Optimization by vector space methods}.\hskip 1em plus
  0.5em minus 0.4em\relax Hoboken, New Jersey: John Wiley \& Sons Inc., 1969.

\bibitem{Cai87}
P.~E. Caines, \emph{Linear stochastic systems}.\hskip 1em plus 0.5em minus
  0.4em\relax New York, NY, USA: John Wiley \& Sons, Inc., 1987.

\bibitem{MarPel08b}
D.~Marinazzo, M.~Pellicoro, and S.~Stramaglia, ``Kernel {G}ranger causality and
  the analysis of dynamical networks,'' \emph{Phys. Rev. E}, vol.~77, p.
  056215, 2008.

\bibitem{GolRee98}
J.~S. Goldstein, I.~S. Reed, and L.~L. Scharf, ``A multistage representation of
  the {W}iener filter based on orthogonal projections,'' \emph{IEEE
  Transactions on Information Theory}, vol.~44, pp. 2943--2959, 1998.

\end{thebibliography}

\end{document}